\documentclass[a4paper]{article} 
\usepackage{amsmath}
\usepackage{amssymb}
\usepackage[margin=4cm]{geometry}
\usepackage{graphicx}
\usepackage{multirow}
\usepackage{url}
\usepackage{color}
\usepackage{comment}
\usepackage{upref}
\usepackage{hyperref}
\usepackage{fancyvrb}
\usepackage{epstopdf}
\usepackage{algorithm} 
\usepackage{algorithmicx}
\usepackage{algpseudocode}
\usepackage{afterpage}
\usepackage{lipsum}

\definecolor{hotpink}{rgb}{0.9,0,0.5}
\hypersetup{colorlinks,urlcolor=blue,citecolor=hotpink,linkcolor=blue}

\def\yh{\widehat{\mathcal{Y}}}
\def\uy{\underline{\mathcal{Y}}}
\def\y{\mathcal{Y}}
\def\K{\mathcal{K}}
\def\P{\mathcal{P}}
\def\R{\mathcal{R}}
\def\lmax{\tilde{\lambda}_{\text{max}}}
\def\lmin{\tilde{\lambda}_{\text{min}}}
\def\xhat{\widehat{x}}
\def\rhat{\widehat{r}}

\def\B{\mathcal{B}}
\def\Gh{\widehat{G}}
\def\veps{\varepsilon}
\def\norm#1{\|#1\|}

\algnewcommand{\IfThen}[2]{
  \State \algorithmicif\ #1\ \algorithmicthen\ #2}

\newcounter{mylineno}
\makeatletter
\let\oldtabcr\@tabcr

\def\mynewline{\refstepcounter{mylineno}%
                \llap{\footnotesize\arabic{mylineno}\hspace{5pt}}%
               }

\gdef\@tabcr{\@stopline \@ifstar{\penalty%
            \@M \@xtabcr}\@xtabcr\mynewline}

\makeatother

\mathchardef\Gamma="7100
\mathchardef\Delta="7101
\mathchardef\Theta="7102
\mathchardef\Lambda="7103
\mathchardef\Xi="7104
\mathchardef\Pi="7105
\mathchardef\Sigma="7106
\mathchardef\Upsilon="7107
\mathchardef\Phi="7108
\mathchardef\Psi="7109
\mathchardef\Omega="710A

\mathcode`@="8000 
{\catcode`\@=\active\gdef@{\mkern1mu}}

\begin{document}

\newtheorem{The}{Theorem}[section]

\numberwithin{equation}{section}
\title{An Adaptive $s$-step Conjugate Gradient Algorithm with Dynamic Basis Updating}

\author{Erin C. Carson\thanks{Faculty of Mathematics and Physics, Charles University. This research was supported by OP RDE project No. CZ.02.2.69/0.0/0.0/16\_027/0008495, International Mobility of Researchers at Charles University.
}}

\date{}

\maketitle
\begin{abstract}
The adaptive $s$-step CG algorithm is a solver for sparse, symmetric positive definite linear systems designed to reduce the synchronization cost per iteration while still achieving a user-specified accuracy requirement. In this work, we improve the adaptive $s$-step conjugate gradient algorithm by use of iteratively updated estimates of the largest and smallest Ritz values, which give approximations of the largest and smallest eigenvalues of $A$, using a technique due to Meurant and Tich{\' y} [G. Meurant and P. Tich{\' y}, Numer. Algs. (2018), pp.~1--32]. The Ritz value estimates are used to dynamically update parameters for constructing Newton or Chebyshev polynomials so that the conditioning of the $s$-step bases can be continuously improved throughout the iterations. These estimates are also used to automatically set a variable related to the ratio of the sizes of the error and residual, which was previously treated as an input parameter. We show through numerical experiments that in many cases the new algorithm improves upon the previous adaptive $s$-step approach both in terms of numerical behavior and reduction in number of synchronizations. 
\end{abstract}




\IfFileExists{gitHeadInfo.gin}%
{\thispagestyle{fancy}
\renewcommand{\headrulewidth}{0pt} 
\fancyhead[LE,LO]{}
\fancyhead[RE,RO]{{
                  \BeforeSubString{+0000}{\gitAuthorIsoDate}}
\gitAbbrevHash}
}{}

%
\section{Introduction}
\label{sec:intro}

In this work we focus on the problem of solving linear systems $Ax=b$, where $A\in \mathbb{R}^{N \times N}$ is symmetric positive definite (SPD). When $A$ is large and sparse, the  iterative conjugate gradient method (CG), which is a Krylov subspace method, is commonly-used as a solver. Given an initial approximate solution $x_0$ with initial residual $r_0 = b-Ax_0$, Krylov subspace methods construct a sequence of nested Krylov subspaces $\K_1(A,r_0) \subset \K_2(A,r_0) \subset \cdots \subset \K_i(A,r_0)$ where 
\[
\K_{i}(A,r_0) = \text{span}\{r_0, Ar_0,\ldots, A^{i-1}r_0\}. 
\]
In CG, the approximate solution $x_{i}\in x_0 + \K_{i}(A,r_0)$ is constructed according to the orthogonality constraint $r_{i} = b- Ax_{i} \perp \K_i(A,r_0)$, which is equivalent to selecting the vector $x_{i} \in x_0 + \K_i(A,r_0)$ that minimizes the $A$-norm (energy norm) of the error, i.e., $\Vert x- x_{i} \Vert_A = \sqrt{(x-x_{i})^T A (x-x_{i})}$.  

Perhaps the most well-known algorithm for CG is due to Hestenes and Stiefel~\cite{hest52}, which uses three coupled 2-term recurrences for recursively updating the approximate solution $x_{i}$, the residual $r_{i}$, and $A$-conjugate search direction vectors $p_{i}$. We refer to this particular algorithm, displayed in Algorithm~\ref{alg:cg}, as ``HSCG'' in this work. 

The CG method is closely related to the Lanczos tridiagonalization method, which iteratively constructs a symmetric tridiagonal matrix $T_i\in \mathbb{R}^{i \times i}$ such that 
$T_i = V_i^T AV_i$ where $V_i$ is an orthonormal basis for the Krylov subspace $\K_i(A,r_0)$. The $i$ eigenvalues of $T_i$, called \emph{Ritz values}, provide estimates of the eigenvalues of $A$. The matrix $T_i$ can be written in terms of the coefficients $\alpha_i$ and $\beta_i$ computed in lines~\ref{alp} and~\ref{bet} of Algorithm~\ref{alg:cg} via
\begin{equation}
T_i = 
\begin{bmatrix}
\frac{1}{\alpha_0} & \frac{\sqrt{\beta_0}}{\alpha_0} & & & \\
\frac{\sqrt{\beta_0}}{\alpha_0} & \frac{1}{\alpha_1} + \frac{\beta_0}{\alpha_0} & \frac{\sqrt{\beta_1}}{\alpha_1} & & \\
& \frac{\sqrt{\beta_1}}{\alpha_1} & \frac{1}{\alpha_2} + \frac{\beta_1}{\alpha_1} & \ddots & \\
& & \ddots & \ddots & \frac{\sqrt{\beta_{i-2}}}{\alpha_{i-2}} \\
& & & \frac{\sqrt{\beta_{i-2}}}{\alpha_{i-2}} & \frac{1}{\alpha_{i-1}} + \frac{\beta_{i-2}}{\alpha_{i-2}}
\end{bmatrix};
\label{eq:T}
\end{equation}
see, e.g.,~\cite[Sec. 6.7.3]{saad03}.

\begin{algorithm}
\caption{Hestenes and Stiefel CG (HSCG)}
\label{alg:cg}
\begin{algorithmic}[1]
\Require {$N \times N$ symmetric positive definite matrix $A$, length-$N$ vector $b$, initial approximation $x_0$ to $Ax=b$, desired convergence criterion} 
\Ensure {Approximate solution $x_{i+1}$ to $Ax=b$ with updated residual $r_{i+1}$}
\State {$r_{0}=b-Ax_{0},\, p_{0}=r_{0}$}
\For {$i=0,1,\dots,$ until convergence}
\State $\alpha_{i}=r^T_i r_i/p^T_i A p_i$ \label{alp}
\State $q_i = \alpha_i p_i$
\State $x_{i+1}=x_i+q_i$
\State $r_{i+1}=r_i-Aq_i$ 
\State $\beta_{i}= r^T_{i+1} r_{i+1}/ r^T_i r_i$ \label{bet}
\State $p_{i+1}=r_{i+1}+\beta_i p_i$ 
\EndFor
\end{algorithmic}
\end{algorithm}

In the setting of large-scale problems on parallel machines, the performance of HSCG is limited by communication, i.e., data movement, due to the sparse matrix vector product (SpMV) and inner products in each iteration, both of which have low computational intensity; see, e.g.,~\cite{dhl16}. This has led to the development of a number of algorithmic variants of CG which aim to reduce the communication and/or synchronization cost over a fixed number of steps. 
One such variant is called $s$-step CG (also called communication-avoiding CG; see, e.g.,~\cite{hoemmen2010communication},~\cite{carsonthesis} and the historical references therein). The $s$-step CG algorithm works by computing $O(s)$ new basis vectors for the Krylov subspace at a time and then computing a block inner product between these computed basis vectors. The former can be accomplished with $O(1)$ messages between parallel processors (usually assumed to involve communication only between neighboring processors) under some constraints on the sparsity structure and parallel partition of the matrix $A$; see~\cite{dhmy08} for details. The latter can be accomplished with a single global synchronization point. Over $s$ iterations, this approach can reduce the number of synchronizations from $O(s)$ to $O(1)$. We elaborate on the mathematics behind this approach in Section~\ref{sec:sstep}.  

In practice, however, we are not concerned only with the number of synchronizations over a fixed number of steps, but rather the total number of synchronizations required to achieved the prescribed convergence criterion. It is well-known that $s$-step CG (and $s$-step variants of other Krylov subspace methods) can exacerbate the delay of convergence and decrease in attainable accuracy that are characteristic of finite precision HSCG. This behavior generally grows worse for increasing $s$. In the extreme case, this can lead to a situation where the prescribed accuracy is no longer attainable, making the $s$-step approach inapplicable. 

In $s$-step CG, it has been shown that the loss of attainable accuracy can be bounded in terms of the condition numbers of the computed $O(s)$-dimensional bases for the Krylov subspaces generated at the beginning of each block of $s$ iterations~\cite{carson2014residual}. This rounding error analysis led to the insight that if we wish to achieve a certain accuracy, then the condition numbers of the bases must be controlled to be less than a certain quantity inversely proportional to the largest residual norm within the current block of $s$ iterations. This naturally suggests a variable $s$-step approach; when the size of the residual is still large (at the beginning of the iterations), $s$ should be small so as to keep the bases well-conditioned, but as the method converges and the size of the residual decreases, $s$ can be gradually increased without detriment to the maximum attainable accuracy. 

This inspired the \emph{adaptive $s$-step CG algorithm}~\cite{carson18}, in which the value of $s$ is automatically adjusted according to the basis condition number and the user-specified accuracy requirement. It was shown in~\cite{carson18} that this adaptive approach can provide improved reliability in terms of numerical behavior while still reducing the overall number of synchronizations. Given the importance of the condition of the bases, one particular gap (which applies to most $s$-step CG algorithms in the literature)
is the reliance on user-supplied parameters for generating the polynomial bases for the $O(s)$-dimensional Krylov subspaces. If this information is not known a priori, either basis parameters must be computed before execution of the $s$-step algorithm (by precomputing information about the spectrum of $A$ or running a number of iterations of HSCG). A simple monomial basis can be used as a default, it is known that the condition number of the monomial basis grows exponentially with $s$~\cite{gautschi79}. 
The previous method also required the user to heuristically set a certain parameter used in determining how large $s$ can be in each outer loop iteration. 

In this work, we develop an improved version of the adaptive $s$-step CG algorithm which makes use of estimates of the largest and smallest eigenvalues of $A$ (the largest and smallest Ritz values) obtained automatically and inexpensively as a byproduct of the iterations using the technique of Meurant and Tich{\' y}~\cite{meti18}. 
The eigenvalue estimates are used for two purposes: (1) to dynamically update coefficients for Newton or Chebyshev polynomials so that the conditioning of the $s$-step bases can be continuously improved throughout the iterations and (2) to automatically set the previously heuristically-chosen parameter based on information obtained during the iterations.

In Section~\ref{sec:sstep}, we briefly review $s$-step CG. We then outline the idea of the adaptive $s$-step CG algorithm in Section~\ref{sec:asstep}. In Section~\ref{sec:iasstep}, we review the work of Meurant and Tich{\' y}~\cite{meti18} and use this to develop an improved adaptive $s$-step CG algorithm. Section~\ref{sec:exp} presents numerical experiments for a variety of small test problems, which demonstrate the benefits of the improved approach. We conclude and discuss further challenges in Section~\ref{sec:conc}. 

\section{The $s$-step CG algorithm}
\label{sec:sstep}

The idea of $s$-step Krylov subspace algorithms is not new. The first known appearance of $s$-step CG in the literature was due to van Rosendale~\cite{rosendale83}, although the name ``$s$-step CG'' was first coined later by Chronopoulos and Gear~\cite{chge89}. There have since been many efforts towards developing $s$-step formulations of CG and other Krylov subspace methods; for a thorough overview of related works, see~\cite[Table 1.1]{hoemmen2010communication}. Much of the early work in this area was motivated by reducing the amount of I/O and/or by increasing the potential for parallelism in CG. These algorithms (and other variants of Krylov subspace methods designed to reduce communication overhead, such as pipelined algorithms~\cite{gamv13}, \cite{ghy14})  have recently come back into vogue for their potential to reduce data movement, both between levels of the memory hierarchy on a single processor and between processors in the parallel setting, which becomes increasingly important in efforts to scale to larger problem sizes and larger machines~\cite{dongarra11}. 

The $s$-step approach can be thought of as blocking the iterations into sets of $s$. For one block of $s$ iterations, one first expands the underlying Krylov subspace by $O(s)$ dimensions and subsequently performs a block orthogonalization using only a single global synchronization point. The vector updates for the block of $s$ iterations can then be performed by updating the $O(s)$ coordinates of the vectors in the generated $O(s)$-dimensional Krylov subspace. The particular details of the algorithm depend on the particular underlying Krylov subspace method, but the general concept is the same. In order to establish notation, we give a brief overview of the $s$-step CG algorithm. 

The $s$-step CG algorithm consists of an outer loop, indexed by $k\geq 0$, which iterates over the blocks of $s$ iterations, 
and an inner loop, which iterates over $j\in\{0,\ldots,s-1\}$ within each block. For clarity, we globally 
index iterations by $i\equiv sk+j$. It follows from the properties of CG that at the beginning of a block $k$, for $\ell\in\{0,\ldots,s\}$ we have
\begin{align}
p_{sk+\ell}, r_{sk+\ell} &\in \K_{\ell+1} (A, p_{sk}) + \K_\ell (A,r_{sk}), \nonumber \\
x_{sk+\ell}-x_{sk} &\in \K_{\ell} (A, p_{sk}) + \K_{\ell-1} (A,r_{sk}). \label{eq:reqsub}
\end{align}
Then the CG vectors for the next $s$ iterations to be computed within this block lie in the union of the column spaces of the matrices
\begin{align}
\P_{k,s} &= [\rho^{(k)}_0(A)p_{sk},\ldots,\rho^{(k)}_s(A)p_{sk}], \quad \text{span}(\P_{k,s}) = \K_{s+1}(A,p_{sk}), \nonumber \\
 \R_{k,s} &= [\rho^{(k)}_0(A)r_{sk},\ldots,\rho^{(k)}_{s-1}(A)r_{sk}], \quad \text{span}(\R_{k,s}) = \K_{s}(A,r_{sk}),
\label{eq:cg-krylovbasis}
\end{align}
where $\rho^{(k)}_\ell(z)$ is a polynomial of degree $\ell$ satisfying the three-term recurrence
\begin{align}
\rho^{(k)}_0(z) &=1, \qquad \rho^{(k)}_1(z) = (z-\theta^{(k)}_0)\rho^{(k)}_0(z)/\gamma^{(k)}_0 \nonumber \\
\rho^{(k)}_\ell(z) &= ((z-\theta^{(k)}_{\ell-1})\rho^{(k)}_{\ell-1}(z) - \mu^{(k)}_{\ell-2}\rho^{(k)}_{\ell-2}(z))/\gamma^{(k)}_{\ell-1}, \quad \ell\geq 2. \label{eq:rho}
\end{align}
The need for the superscripts $(k)$ above will become clear later when we introduce the improved adaptive $s$-step CG algorithm, in which the coefficients in the recurrence~\eqref{eq:rho} are dynamically updated between outer loop iterations. 
Under certain constraints on the sparsity structure and the partition of the data, the matrices~\eqref{eq:cg-krylovbasis} can be computed with $O(1)$ messages per processor in a parallel algorithm; see~\cite{dhmy08} for details.

We define the \emph{$s$-step basis matrix} $\y_{k,s} = [\P_{k,s},\R_{k,s}]$ and we define $\uy_{k,s}$ to be the same as $\y_{k,s}$ except with columns 
$s+1$ and $2s+1$ set to zero. We can then write the recurrence relation 
\begin{equation}
A \uy_{k,s} = \y_{k,s} \B_{k,s}, 
\label{eq:AVVB}
\end{equation}
where
\begin{equation}
\B_{k,s} = 
\begin{bmatrix}
B^{(k)}_{s+1} & 0 \\
0 & B^{(k)}_s
%
\end{bmatrix}, 
\quad\text{with} \quad B^{(k)}_{i} \equiv 
	\begin{bmatrix}
	\theta^{(k)}_0 & \mu^{(k)}_0 & & & \\
	\gamma^{(k)}_0 & \theta^{(k)}_1 & \ddots & & \\
	& \gamma^{(k)}_1 & \ddots & \mu^{(k)}_{i-3} & \\
	& & \ddots & \theta^{(k)}_{i-2} & 0 \\
	& & & \gamma^{(k)}_{i-2} & 0 
	\end{bmatrix}.
\label{B}
\end{equation}

For $\ell\in\{0,\ldots, s\}$ we can then represent the vectors $x_{sk+\ell}-x_{sk}$ , $r_{sk+\ell}$, and $p_{sk+\ell}$ by their $2s+1$ coordinates in the 
basis spanned by the columns of $\y_{k,s}$, i.e.,
\begin{equation}
[ x_{sk+\ell}-x_{sk}, r_{sk+\ell}, p_{sk+\ell}] = \y_{k,s}[x'_{k,\ell},r'_{k,\ell}, p'_{k,\ell} ],
\label{eq:coords} \\
\end{equation}
and the updates to these coordinate vectors in the inner loop from $j\in\{0,\ldots,s-1\}$ become
\begin{align*}
x'_{k,j+1} &= x'_{k,j} + \alpha_{sk+j}p'_{k,j}, \\
r'_{k,j+1} &= r'_{k,j}-\alpha_{sk+j}\B_{k,s}p'_{k,j}, \\
p'_{k,j+1} &= r'_{k,j+1} + \beta_{sk+j} p'_{k,j}.
\end{align*}
In practical applications we expect $s\ll N$, so the updates to the length-$(2s+1)$ coordinate vectors can be accomplished locally on each processor without any further communication. 
The coefficients $\alpha_{sk+j}$ and $\beta_{sk+j}$ can also be computed locally on each processor without communication; using~\eqref{eq:AVVB} and~\eqref{eq:coords}, we have  
\begin{align*}
\alpha_{sk+j} &= \frac{r_{sk+j}^T r_{sk+j}}{p_{sk+j}^T A p_{sk+j}} = \frac{(\y_{k,s} r'_{k,j})^T (\y_{k,s} r'_{k,j})}{(\y_{k,s} p'_{k,j})^T (\y_{k,s}\B_{k,s} p'_{k,j})}
= \frac{r'^T_{k,j} (\y_{k,s}^T \y_{k,s}) r'_{k,j}}{p'^T_{k,j}  (\y_{k,s}^T \y_{k,s})\B_{k,s} p'_{k,j}} \\&= \frac{r'^T_{k,j}G_{k,s}r'_{k,j}}{p'^T_{k,j} G_{k,s}\B_{k,s}p'_{k,j}}, \\
\beta_{sk+j} &= \frac{r_{sk+j+1}^T r_{sk+j+1}}{r_{sk+j}^T r_{sk+j}} = \frac{(\y_{k,s} r'_{k,j+1})^T (\y_{k,s} r'_{k,j+1})}{(\y_{k,s} r'_{k,j})^T (\y_{k,s} r'_{k,j})}
= \frac{r'^T_{k,j+1} (\y_{k,s}^T \y_{k,s}) r'_{k,j+1}}{r'^T_{k,j} (\y_{k,s}^T \y_{k,s}) r'_{k,j}} \\ &= 
\frac{r'^T_{k,j+1} G_{k,s}r'_{k,j+1}}{r'^T_{k,j}G_{k,s}r'_{k,j}},
\end{align*}
where 
\[
G_{k,s} = \y_{k,s}^T \y_{k,s}
\]
is the $(2s+1)\times(2s+1)$ Gram matrix which is computed only once per outer loop iteration (requiring a single global synchronization) and stored locally on each processor. 

In~\eqref{B}, the choice $\theta^{(k)}_\ell=0$, $\gamma^{(k)}_\ell=1$, and $\mu^{(k)}_\ell=0$ for all $\ell$ corresponds to the monomial basis. Without 
further information about the spectrum, this is the simplest choice we have, although it is known that the monomial 
basis quickly becomes ill-conditioned with $s$~\cite{gautschi79}, and this thus limits the $s$ we can choose.
Already in early works on $s$-step CG it was observed that the conditioning of the $s$-step basis matrices plays a large role in the resulting finite precision behavior, which led many to experiment with more well-conditioned polynomial bases such as Newton or Chebyshev bases; see, e.g,~\cite{hiwa86, de91, joca92, bhr94, erhel95, deva95}. 

 Estimates of the maximum and minimum eigenvalues of $A$, $\lmax$ and $\lmin$, 
can be used to construct either Newton or Chebyshev polynomials, which will in general result in a better-conditioned basis and thus allow the use of larger $s$ values without loss of accuracy. We note that the idea of adaptively improving basis conditioning using Ritz values is not new. 
A method for generating parameters for Newton and Chebyshev polynomial bases for Krylov subspaces based on Ritz values is described in~\cite{phre12}; see also the works~\cite{manteuffel78} and~\cite{cgr94}. The modern software package Trilinos~\cite{trilinos} currently includes the ability to automatically generate Newton or Chebyshev polynomials in their 
implementations of (fixed) $s$-step Krylov subspace algorithms. The approach they use is to first perform a fixed number of iterations of the classical algorithm (with $s=1$) and then  
use the Ritz value estimates (eigenvalues of $T_i$) to generate basis parameters.

We note that, in contrast with the use of estimates of $\lmax$ and $\lmin$ within the Chebyshev semi-iterative method, for the purposes of generating a Krylov subspace basis, $\lmax$ and $\lmin$ do not need to be particularly accurate in order to provide an improvement over the monomial basis; see, e.g.,~\cite[pp.~12]{phre12}, where it is stated (referring to the Chebyshev basis) that ``it is not important that the ellipse be determined to high accuracy.'' 

Given $\lmax$ and $\lmin$, the Newton basis parameters can be taken as (for now dropping the superscripts) $\mu_\ell = 0$, $\gamma_\ell=1$, and  
\begin{align}
\theta_0 &= \lmax, \qquad \theta_1 = \lmin, \nonumber \\
\theta_\ell &= \text{argmax}_{\theta\in [\lmin,\lmax]} \prod_{m=0}^{\ell-1} |\theta-\theta_m|, \qquad \text{for}\quad \ell\in\{2,\ldots,s-1\}, \label{eq:newtbasis}
\end{align}
which corresponds to a Leja ordering of the points on the real line between $\lmin$ and $\lmax$. The Leja ordering is known to improve the accuracy of operations on polynomials; see, e.g.,~\cite{bhr94},~\cite{care03}, and~\cite{reichel90}. 
Parameters for a simplified version of the Chebyshev basis can be chosen as (see \cite[Section 4.4]{joca92})
\begin{align}
\theta_\ell &= \frac{\lmin+\lmax}{2} \qquad\text{for}\quad \ell\in\{0,\ldots,s-1\}, \nonumber \\
\qquad \mu_\ell &= 2(\lmax-\lmin) \qquad\text{for}\quad \ell\in\{0,\ldots,s-2\}, \quad\text{and} \nonumber \\
\gamma_0 &= \lmax-\lmin, \quad \gamma_\ell = \frac{\lmax-\lmin}{2} \qquad\text{for}\quad \ell\in\{1,\ldots,s-1\}. \label{eq:chebbasis}
\end{align}

\subsection{The $s$-step CG algorithm in finite precision}

It is well-known that finite precision roundoff errors can cause a delay of convergence and a decrease of attainable accuracy in CG algorithms; see, e.g.,~\cite{mest06},~\cite{list12} and references therein. 
In this work, we focus on the maximum attainable accuracy, as achieving a certain prescribed accuracy is the goal of the adaptive $s$-step approaches. 

The size of the true residual is often used as a computable measure of accuracy. The mechanism by which accuracy is lost is the deviation of the recursively updated residual $\rhat_{i}$ and the true residual $b-A\xhat_{i}$, where $\rhat_{i}$ and $\xhat_{i}$ denote the quantities computed in finite precision (in general, we will now use hats to denote quantities computed in finite precision). Writing $b-A\xhat_{i} =  (b-A\xhat_{i} -\rhat_{i}) + \rhat_{i}$, it is clear that as the size of the recursively updated residual becomes very small, the upper bound on the size of $b-A\xhat_{i}$ depends on the size of the \emph{residual gap} $\delta_{i} \equiv b-A\xhat_{i} -\rhat_{i}$. 

There is a large literature on analyses of maximum attainable accuracy in HSCG and related CG algorithms, including the works of Greenbaum~\cite{greenbaum97}, van der Vorst and Ye~\cite{vaye99}, Sleijpen and van der Vorst~\cite{slvo96}, and Gutknecht and Strako{\v s}~\cite{gust00}. 
Modifying slightly the bound on the residual gap in HSCG derived by Sleijpen and van der Vorst~\cite[Eqn. (6)]{slvo96}, we can bound the growth of the residual gap starting from some iteration $m\equiv sk$ to iteration $m+j+1$, $j\in\{0,\ldots,s-1\}$, in HSCG by
\begin{equation}
\Vert \delta_{m+j+1} -\delta_{m} \Vert \leq c \veps  \left(\max_{0\leq \ell \leq j+1} \Vert \rhat_{m+\ell}\Vert \right),
\label{eq:cggap}
\end{equation} 
where $c= 2sN_A \nu \kappa(A)$, $\kappa(A) = \Vert A^{-1} \Vert \Vert A \Vert$ denotes the condition number, $N_A$ is the maximum number of nonzeros per row in $A$, $\nu =  \Vert \vert A \vert \Vert / \Vert A \Vert$, $\Vert \cdot \Vert$ denotes the 2-norm, and $\veps$ is the machine unit roundoff. 

In~\cite[Eqn. (16)]{carson18} it is shown that the growth of the residual gap within one outer loop iteration $k$ (which begins at global iteration $m$) of $s$-step CG can be bounded by 
\begin{equation}
\Vert \delta_{m+j+1} - \delta_{m}\Vert \leq c_{k} \varepsilon \kappa(\yh_{k,s}) \left( \max_{0\leq \ell \leq j+1} \Vert \rhat_{m+\ell}\Vert \right) + \varepsilon \Vert A \Vert \Vert x \Vert,
\label{eq:growth}
\end{equation}
for $j\in\{0,\ldots,s-1\}$. 
The constant $c_{k}$ comes from the rounding error analysis and can be written as 
\begin{equation}
c_{k} = 2s \big( 2(3+N_A)\nu t + (6+8t)\tau_k + 2t^3 + 3 \big) \kappa(A),
\label{eq:ck}
\end{equation}
where $t = \sqrt{2s+1}$ and $\tau_k = \Vert \vert \B_{k,s} \vert \Vert / \Vert A \Vert$. 
The most notable difference between~\eqref{eq:cggap} and~\eqref{eq:growth} is the appearance of the term $\kappa(\yh_{k,s})=\Vert \yh_{k,s}^{+} \Vert \Vert \yh_{k,s} \Vert$, where $\yh_{k,s}^{+}$ denotes the Moore-Penrose pseudoinverse. 
In other words, the local roundoff errors made in $s$-step CG are \emph{amplified} by the condition numbers of the computed $s$-step basis matrices. 
This theoretically confirms observations regarding the effect of the condition numbers of the $s$-step bases on the numerical behavior of $s$-step CG compared to HSCG; an ill-conditioned $s$-step basis can cause an increase in the residual gap and thus can decrease the attainable accuracy.

\section{The adaptive $s$-step CG algorithm}
\label{sec:asstep}

Rearranging \eqref{eq:growth}, it can be shown that if the application requires a relative residual norm of $\veps^*$, the condition number of the basis matrix $\yh_{k,s}$ must satisfy
\begin{equation}
\kappa(\yh_{k,s}) \leq \frac{\veps^*}{c_{k} \veps \Vert \rhat_{m+\ell}\Vert}, \qquad \text{for}\quad 0 \leq \ell \leq s.
\label{cond}
\end{equation}
This naturally suggests that $s$ should be allowed to vary in each outer loop $k$; when the residuals are large, $\kappa(\yh_{k,s})$ and thus $s$ must be small, but as the residual is reduced, the condition number of the basis and thus $s$ can be larger without detriment to the attainable accuracy\footnote{
We note that this insight is similar to that behind the development of the so-called ``inexact Krylov subspace methods''; see, e.g.,~\cite{ss03} as well as the technical report~\cite{bv00} which was later published as~\cite{bf05}. 
}. We therefore introduce the subscripted quantity $s_k$ to denote the number of inner iterations in outer loop $k$. Quantities which depend on $s_k$ will now also have a subscript $k$, e.g., $\yh_{k,s_k}$, $\B_{k,s_k}$, $t_k = \sqrt{2s_k+1}$, and the expression for $c_k$ now containing the relevant quantities dependent on the value of $s_k$ in outer loop $k$.  

This idea led to the adaptive $s$-step CG algorithm published in~\cite{carson18}, displayed in Algorithm~\ref{alg:vscg}.  
It is shown in~\cite{carson18} that, assuming that the algorithm converges, as long as \eqref{cond} is satisfied, then adaptive $s$-step CG can attain a solution to the same level of accuracy as HSCG. We now give a brief description of the algorithm.


In CG there is no guarantee that the residual norms are monotonically decreasing (in exact arithmetic, the method rather minimizes the $A$-norm of the error). 
Therefore, for an outer loop beginning at iteration $m=\sum_{\ell=0}^{k-1} s_\ell$, we use $\bar{s}_k$ to denote our initial guess for $s_k$ and construct the basis $\yh$ only based on the current residual $\rhat_{m}$. We then compute the Gram matrix $\Gh_{k,\bar{s}_k}=\yh_{k,\bar{s}_k}^T \yh_{k,\bar{s}_k}$. From this, we find the largest value $\tilde{s}_k \leq \bar{s}_k$ such that condition~\eqref{cond} holds, i.e., $\tilde{s}_k$ is the maximum value in $\{1,\ldots,\bar{s}_k\}$ such that 
 \begin{equation}
\kappa(\yh_{k,\tilde{s}_k})  \leq \frac{\veps^*}{c_{k}\veps\norm{\rhat_{m}}} 
\label{condy}
\end{equation}
holds. The values $\kappa(\yh_{k,i})$ can be estimated using the square roots of the leading principle submatrices of the constructed $\Gh_{k,\bar{s}_k}$, since 
since $\kappa(\yh_{k,\ell}) \approx \sqrt{\kappa(\Gh_{k, \ell})}$. The quantities $\sqrt{\kappa(\Gh_{k, \ell})}$ are inexpensive to compute; this involves $O(s^3)$ floating point operations to compute the eigenvalues of $\Gh_{k, \bar{s}_k}$ and no additional data movement as $\Gh_{k, \bar{s}_k}$ is stored locally on each processor. 
We then check if the condition~\eqref{condy} is violated within each inner loop, which can occur if we encounter a large intermediate residual norm. We use $s_k \leq \tilde{s}_k$ to denote the \emph{actual} number of inner loop iterations which occurred. 
We note that the residual norms can be estimated cheaply (without communication) within the inner loop since (in exact arithmetic)
\[
\Vert r_{m+j+1} \Vert = \sqrt{r'^T_{k,j+1} G_{k,\tilde{s}_k} r'_{k,j+1}}.
\]

The algorithm also requires the user to input some $\bar{s}_0$ as an initial value, a value $\sigma$ which is the maximum value for $\bar{s}_k$ (which could be determined by offline auto-tuning and should be based on the matrix nonzero structure, machine parameters, and matrix partition), and a value $f$, which is the maximum $\bar{s}_k$ can be allowed to grow in each iteration (e.g., $f$ could be made small to reduce the amount of wasted flops if $s_k$ is much smaller than $\bar{s}_k$ in each outer loop, or $f$ could be made larger to maximize the potential $s_k$ in each outer loop). 

In the experiments in~\cite{carson18}, it was found that the value of $c_k$ in~\eqref{eq:ck} is often a large overestimate, resulting in smaller $s_k$ values than necessary to achieve the desired accuracy which results in more outer loop iterations than necessary. It was found that in most cases, taking $c_k=1$ in~\eqref{cond} worked well, although there was no theoretical justification to support this. In the following section, we describe a way to adaptively and automatically set this parameter based on existing quantities obtained during the iterations.

We briefly comment on related work in the area of using a variable $s$ value in $s$-step Krylov subspace algorithms. Also motivated by improving numerical behavior, Imberti and Erhel used a variable $s$ value in their $s$-step GMRES algorithm~\cite{imer17}, although their approach requires the user to prescribe a priori the sequence of $s_k$ values. A variable $s$ value was also used within a $s$-step BICGSTAB algorithm used as the coarse grid solve routine within a geometric multigrid method~\cite{williams14}. In~\cite{williams14} this approach was termed a ``telescoping $s$'', in which the value of $s$ starts small and is allowed to grow as the outer loops proceed. This was done for performance reasons rather than numerical ones; when the coarse grid problem is easy (converges in a few iterations), we do not waste effort computing a larger $s$-step basis than need be.


\begin{algorithm}
\caption{Adaptive $s$-step conjugate gradient}
\label{alg:vscg}
\begin{algorithmic}[1]
\Require {$N \times N$ symmetric positive definite matrix $A$, length-$N$ vector $b$, initial approximation $x_0$ to $Ax=b$, maximum $s_k$ value $\sigma$, initial value $\bar{s}_0$, maximum basis growth factor $f$, desired convergence tolerance $\varepsilon^*$, function $c_{k}$} 
\Ensure {Approximate solution $x_{m}$ to $Ax=b$ with updated residual $r_{m}$}
\State {$r_{0}=b-Ax_{0},\, p_{0}=r_{0}$, $m = 0$}
\For {$k=0,1,\dots,$ until convergence}
\IfThen {$k\neq 0$}{$\bar{s}_k = \min(s_{k-1}+f, \sigma)$}
\State {Compute $\bar{s}_k$-step basis matrix $\mathcal{Y}_{k,\bar{s}_k}=[\mathcal{P}_{k,\bar{s}_k},\, \mathcal{R}_{k,\bar{s}_k}]$ according to~\eqref{eq:cg-krylovbasis}.}
\State {Compute $G_{k,\bar{s}_k}=\mathcal{Y}_{k,\bar{s}_k}^T\mathcal{Y}_{k,\bar{s}_k}$.} 
\State {Determine $\tilde{s}_k$ by~\eqref{condy}; assemble $\mathcal{Y}_{k,\tilde{s}_k}$ and $G_{k,\tilde{s}_k}$.}
\State {Store estimate $\gamma \approx \kappa(\mathcal{Y}_{k,\tilde{s}_k})$.}
\State Assemble $\mathcal{B}_{k,\tilde{s}_k}$ such that~\eqref{eq:AVVB} holds.

\State $p'_{k,0}=\left[1,\,0_{1,2\tilde{s}_k}\right]^{T}$, $r'_{k,0}=\left[0_{1,\tilde{s}_k+1},\,1,\,0_{1,\tilde{s}_k-1}\right]^{T}$, $x'_{k,0}=\left[0_{1,2\tilde{s}_k+1}\right]^T$
\For{$j = 0$ to $\tilde{s}_k-1$}
\State $s_k=j+1$
\State $\alpha_{m+j}=\big({r}'^T_{k,j}G_{k,\tilde{s}_k}r'_{k,j}\big)/\big({p}'^T_{k,j}G_{k,\tilde{s}_k}\mathcal{B}_{k,\tilde{s}_k}p'_{k,j}\big)$ 
\State $q'_{k,j} = \alpha_{m+j}p'_{k,j}$

\State $x'_{k,j+1}=x'_{k,j}+q'_{k,j}$
\State $r'_{k,j+1}=r'_{k,j}-{\mathcal{B}_{k,\tilde{s}_k}}q'_{k,j}$ 

\State $\beta_{m+j}=\big({r}'^T_{k,j+1}G_{k,\tilde{s}_k}r'_{k,j+1}\big)/\big({r}'^T_{k,j}G_{k,\tilde{s}_k}r'_{k,j}\big)$ 

\State $p'_{k,j+1}=r'_{k,j+1}+\beta_{m+j}p'_{k,j}$ 

\IfThen {$\gamma \geq \frac{\veps^*}{ c_{k} \veps \left(r_{k,j+1}^{'T} G_{k,\tilde{s}_k} r'_{k,j+1} \right)^{1/2}}$}{break from inner loop.}

\EndFor

\State {Recover iterates $\{p_{m+s_k},r_{m+s_k},x_{m+s_k}\}$ according to~\eqref{eq:coords}.}
\State {$m=m+s_k$}

\EndFor

\end{algorithmic}
\end{algorithm}

\section{An improved adaptive $s$-step CG algorithm}
\label{sec:iasstep}

As stated, the adaptive $s$-step CG algorithm described in~\cite{carson18} requires the user to supply parameters for the polynomial recurrence~\eqref{eq:rho} used in constructing the $s$-step basis matrices as well as the function $c_k$ used in determining when and if to break from the inner loop. In this section we present an improved adaptive $s$-step CG algorithm which uses Ritz value estimates computed via the results of Meurant and Tich{\' y}~\cite{meti18}, which we now briefly summarize. The $\alpha$ and $\beta$ coefficients computed during the CG iterations can be composed to form the Cholesky factor $L_i^T$ of the Lanczos tridiagonal matrix (see~\eqref{eq:T}) $T_{i} = L_i L_i^T$, where 
\begin{equation}
L_i^T \equiv 
\begin{bmatrix}
\zeta_0 & \eta_0 & & \\
 & \ddots & \ddots & \\
& & \ddots & \eta_{i-2}\\
& & & \zeta_{i-1} 
\end{bmatrix}
= 
\begin{bmatrix}
\frac{1}{\sqrt{\alpha_0}} & \sqrt{\frac{\beta_0}{\alpha_0}} & & \\
 & \ddots & \ddots & \\
& & \ddots & \sqrt{\frac{\beta_{i-2}}{\alpha_{i-2}}}\\
& & & \frac{1}{\sqrt{\alpha_{i-1}}}
\end{bmatrix},
\label{eq:lit}
\end{equation} 
with $i$ now denoting the global iteration index (the total number of inner loop iterations). It is expected that the eigenvalues of $T_i$, i.e., the Ritz values, give decent approximations for the extremal eigenvalues of $A$. The extremal eigenvalues of $A$ can be estimated via the relations
\begin{equation}
\lmax = \lambda_{\text{max}}(T_i) = \Vert L_i \Vert^2, \qquad \lmin = \lambda_{\text{min}}(T_i)=\Vert L_i^{-1} \Vert^{-2}.
\label{eq:lest}
\end{equation}
As described in~\cite{meti18}, the norms $\Vert L_i \Vert^2$ and $\Vert L_i^{-1} \Vert^{-2}$ can computed using incremental norm estimation of the matrices $L_i^T$ and $L_i^{-T}$ without the need to explicitly construct these matrices; this requires  
insignificant extra work and no communication. The algorithms for estimating $\Vert L_i \Vert^2$ and $\Vert L_i^{-1} \Vert^{-2}$ are displayed as Algorithms~\ref{alg:L} and~\ref{alg:Li}, respectively, which appear as Algorithms 4 and 5 in~\cite[Sec. 5]{meti18} with slightly modified notation and indexing. 

Thus in each inner loop iteration, we can update estimates of the extremal eigenvalues of $A$, $\lmin$ and $\lmax$, using Algorithms~\ref{alg:L} and~\ref{alg:Li}. We note that in each inner loop iteration, updating $\lmin$ and $\lmax$ requires executing only a single additional for-loop in Algorithms~\ref{alg:L} and~\ref{alg:Li} and storing only a small number of scalar quantities from the previous iteration. 

The estimates $\lmin$ and $\lmax$ are then used in two ways to improve the adaptive $s$-step CG algorithm. The first improvement comes from improving the quality of the computed Krylov subspace bases. The computation begins by using a monomial basis to construct the $\bar{s}_k$-step polynomial bases. Once at least two iterations have finished (i.e., $\lmin \neq \lmax$), the estimates $\lmin$ and $\lmax$ can be used in the subsequent outer loop $k$ to estimate the basis parameters $\theta_i^{(k)}$, $\gamma_i^{(k)}$, and $\mu_i^{(k)}$ for either the Newton basis~\eqref{eq:newtbasis} or the Chebyshev basis~\eqref{eq:chebbasis}. The process continues, with the latest updated estimates $\lmin$ and $\lmax$ being used to generate the basis parameters after every outer loop, improving the quality of the bases as the iterations proceed.

\begin{algorithm}
\caption{Incremental estimation of $\Vert L_i \Vert^2$ (\cite[Alg. 4]{meti18})}
\label{alg:L}
\begin{algorithmic}[1]
\Require {Entries $\{\zeta_0,\ldots, \zeta_{i-1}\}$ and $\{\eta_0,\ldots, \eta_{i-2}\}$ of upper bidiagonal matrix $L_i^T$ in~\eqref{eq:lit}.} 
\Ensure {Quantity $\omega^{\text{max}}_{i-1}$ which gives an estimate of $\Vert L_i \Vert^2$}
\State {$\omega_0 = \zeta_0^2$, $\omega^{\text{max}}_{0}=\omega_0$, $h_0 = 1$}
\For {$\ell = 0,\ldots,{i-2}$}
\State{ $d_\ell = \zeta_\ell^2 \eta_\ell^2 h_\ell$, $a_\ell = \eta_\ell^2 + \xi_{\ell+1}^2$}
\State {$\chi_\ell = \sqrt{(\omega_\ell - a_\ell)^2 + 4 d_\ell}$}
\State {$h_{\ell+1} = \frac{1}{2}\left(1 - \frac{\omega_\ell - a_\ell}{\chi_\ell} \right)$}
\State {$\omega_{\ell+1} = \omega_\ell + \chi_\ell h_{\ell+1}$}
\State {$\omega^{\text{max}}_{\ell+1} = \omega_{\ell+1}$}
\EndFor
\end{algorithmic}
\end{algorithm}

\begin{algorithm}
\caption{Incremental estimation of $\Vert L_i^{-1} \Vert^{-2}$ (\cite[Alg. 5]{meti18})}
\label{alg:Li}
\begin{algorithmic}[1]
\Require {Entries $\{\zeta_0,\ldots, \zeta_{i-1}\}$ and $\{\eta_0,\ldots, \eta_{i-2}\}$ of upper bidiagonal matrix $L_i^T$ in~\eqref{eq:lit}.} 
\Ensure {Quantity $\omega^{\text{min}}_{i-1}$ which gives an estimate of $\Vert L_i \Vert^{-2}$}
\State {$\omega_0 = \zeta_0^{-2}$, $\omega^{\text{min}}_{0}=\zeta_0^2$, $a_0 = \omega_0$, $d_0 = 0$, $g_0 = 0$, $h_0=1$}
\For {$\ell = 0,\ldots,{i-2}$}
\State{ $d_{\ell+1} = -\frac{\eta_\ell}{\zeta_{\ell+1}}(g_\ell d_\ell + h_\ell a_\ell)$}
\State {$a_{\ell+1} = \frac{1}{\zeta_{\ell+1}^2} (\eta_{\ell}^2 a_\ell +1) $}
\State {$\chi_\ell = \sqrt{(\omega_\ell - a_{\ell+1})^2 + 4 d_{\ell+1}^2 } $}
\State {$h_{\ell+1} = \sqrt{\frac{1}{2} \left(  1 - \frac{\omega_\ell - a_{\ell+1} }{\chi_{\ell}} \right)} $}
\State {$\omega_{\ell+1} = \omega_\ell + \chi_\ell h_{\ell+1}^2$}
\State {$g_{\ell+1} = \sqrt{1-h_{\ell+1}^2}$, $h_{\ell+1} = \vert h_{\ell+1} \vert \text{sign}(d_{\ell+1})$}
\State {$\omega^{\text{min}}_{\ell+1} = \omega_{\ell+1}^{-1}$}
\EndFor
\end{algorithmic}
\end{algorithm}

The second improvement comes from using the estimates $\lmin$ and $\lmax$ to eliminate the heuristic choice of the parameter $c_k$.  
As mentioned, the adaptive $s$-step CG algorithm in~\cite{carson18} required that the quantity $c_k$ in the condition~\eqref{cond} be set by the user, as using the full value from~\eqref{eq:ck} was in most cases too restrictive. This is because the $\kappa(A)$ term arises in the quantity $c_k$ as defined in~\eqref{eq:ck} as a result of the effort to bound the residual gap solely in terms of the size of the residuals, for example, 
\begin{align}
\Vert A \Vert \Vert \xhat_{m+j+1} - x\Vert &= \Vert A \Vert \Vert A^{-1}A (\xhat_{m+j+1} - x)\Vert \nonumber \\
&\leq \kappa(A) \Vert A(\xhat_{m+j+1} - x)\Vert = \kappa(A) \Vert \rhat_{m+j+1} \Vert +O(\varepsilon),
\label{eq:b1}
\end{align}
which can give a rather pessimistic bound in practice. We note that a similar observation was made by Sleijpen and van der Vorst~\cite{slvo96}, who argued that the $\kappa(A)$ term will only appear in `unusual' cases. The numerical experiments in~\cite{carson18} found that using $c_k=1$ (i.e., ignoring the $\kappa(A)$ term) worked well in most (but not all) cases.  

We make further use of the results of Meurant and Tich{\' y}~\cite{meti18} and develop a way to adaptively set $c_k$ that removes the burden of setting this parameter heuristically. Instead of the potentially loose upper bound~\eqref{eq:b1}, we will define $\xi_{m+j+1}$ to be the exact quantity defined by 
$\Vert A \Vert \Vert \xhat_{m+j+1} - x\Vert  = \xi_{m+j+1} \Vert \rhat_{m+j+1} \Vert$, i.e., 
\begin{equation}
\xi_{m+j+1} \equiv \frac{\Vert A \Vert \Vert \xhat_{m+j+1} - x\Vert}{\Vert \rhat_{m+j+1} \Vert},
\label{eq:xi}
\end{equation}
where we will have $1\leq \xi_{m+j+1} \leq \kappa(A)$. It is shown in~\cite[Section 3]{meti18} that the $A$-norm of the error can be upper bounded as 
\[
\Vert \xhat_{m+j+1} - x\Vert \leq \frac{\Vert \rhat_{m+j+1} \Vert }{\mu^{1/2}} \psi^{1/2}_{m+j+1},
\]
where $0<\mu\leq \lambda_{\text{min}}$ and the quantity $\psi_{m+j+1} \equiv \Vert \rhat_{m+j+1} \Vert^2 / \Vert \widehat{p}_{m+j+1}\Vert^2$ can be incrementally updated in each iteration by 
\begin{equation}
\psi_{m+j+1} = \frac{\psi_{m+j}}{\psi_{m+j}+\beta_{m+j}}, \quad \psi_0 = 1.
\label{eq:psiupd}
\end{equation}
Thus in each iteration we take $\mu = \lmin$ as computed in~\eqref{eq:lest} and approximately bound $\xi_{m+j+1}$ in~\eqref{eq:xi} via 
\begin{equation}
\tilde{\xi}_{m+j+1} \lesssim \lmax \sqrt{\frac{\psi_{m+j+1}}{\lmin}}.
\label{eq:xiupd}
\end{equation}
Thus in each iteration we set $c_{m+j+1}=\tilde{\xi}_{m+j+1}$ (note the subscript notation change which indicates that this quantity can now change in each inner iteration). 
Updating $\tilde{\xi}_{m+j+1}$ in~\eqref{eq:xiupd} requires only the additional scalar operation in~\eqref{eq:psiupd} to update $\psi_{m+j+1}$. 
Our experiments in Section~\ref{sec:exp} confirm that $\tilde{\xi}_{m+j+1}$ can be far below $\kappa(A)$.

We make a fur\-ther small im\-prove\-ment to the adap\-tive $s$-step CG al\-gor\-ithm in~\cite{carson18}. If we are at the end of some inner iteration $j<\tilde{s}_k-1$, notice that by~\eqref{eq:reqsub}, the iterates that will be updated in the next iteration $j+1$ depend on only a subset of the basis vectors. 
The previous adaptive $s$-step CG algorithm (Algorithm~\ref{alg:vscg}) breaks from the inner loop if in some inner iteration $j$, 
\[
\kappa(\mathcal{Y}_{k,\tilde{s}_k})\geq \frac{\veps^*}{c_{k} \veps \left(r_{k,j+1}^{'T} G_{k,\tilde{s}_k} r'_{k,j+1} \right)^{1/2}},
\]
in other words, if the residual norm $\Vert r_{m+j+1}\Vert$ is so large that~\eqref{cond} will not hold in \emph{some} future inner iteration (not only just the next one). This is overly pessimistic and can cause the algorithm to quit the inner loop iterations unnecessarily early. We can easily and inexpensively modify this approach as follows. 

In each inner loop iteration, we keep track of the maximum residual norm that we have encountered so far in outer loop $k$, and store this as a variable $\phi$. At the end of inner loop iteration $j<\tilde{s}_k-1$, we check whether 
\begin{equation}
\kappa(\mathcal{Y}_{k,j+1})\geq \frac{\veps^*}{c_{m+j+1} \veps \phi},
\label{eq:newcondy}
\end{equation}
and if so, we break from the current inner loop. If not, we continue with the next inner loop iteration $j+1$. 
As previously described, it is easy and inexpensive (involving no communication) to estimate all the condition numbers $\kappa(\yh_{k,\ell+1})$ for $\ell\in\{1,\ldots, \tilde{s}_k-1\}$ by computing the square roots of the condition numbers of the appropriate leading principal submatrices of $\Gh_{k,\tilde{s}_k}$. 
The resulting improved adaptive $s$-step CG algorithm is displayed 
as Algorithm~\ref{alg:vscg2}. The changes compared with Algorithm~\ref{alg:vscg} are highlighted in red.

\begin{algorithm}
\caption{Improved adaptive $s$-step conjugate gradient}
\label{alg:vscg2}
\begin{algorithmic}[1]
\Require {$N \times N$ symmetric positive definite matrix $A$, length-$N$ vector $b$, initial approximation $x_1$ to $Ax=b$, maximum $s$ value $\sigma$, initial $s$ value $\bar{s}_0$, maximum basis growth factor $f$, desired convergence tolerance $\varepsilon^*$} 
\Ensure {Approximate solution $x_{m}$ to $Ax=b$ with updated residual $r_{m}$}
\State {$r_{0}=b-Ax_{0},\, p_{0}=r_{0}$, {\color{red}$\psi_0 = 1$}, {\color{red}$c_1 = \varepsilon^{-1/2}$}, $m = 0$ }
\For {$k=0,1,\dots,$ until convergence}
\IfThen {$k\neq 0$}{$\bar{s}_k = \min(s_{k-1}+f, \sigma)$}
\State {Compute $\bar{s}_k$-step basis matrix $\mathcal{Y}_{k,\bar{s}_k}=[\mathcal{P}_{k,\bar{s}_k},\, \mathcal{R}_{k,\bar{s}_k}]$ according to~\eqref{eq:cg-krylovbasis}.}
\State {Compute $G_{k,\bar{s}_k}=\mathcal{Y}_{k,\bar{s}_k}^T\mathcal{Y}_{k,\bar{s}_k}$.} 
\State {Determine $\tilde{s}_k$ by~\eqref{condy}; assemble $\mathcal{Y}_{k,\tilde{s}_k}$ and $G_{k,\tilde{s}_k}$.}
\State {{\color{red}Store estimates $\gamma_\ell = \kappa(\yh_{k,\ell+1})$ for $\ell\in\{1,\ldots, \tilde{s}_k-1 \}$.}}
\State Assemble $\mathcal{B}_{k,\tilde{s}_k}$ such that~\eqref{eq:AVVB} holds.

\State $p'_{k,0}=\left[1,\,0_{1,2\tilde{s}_k}\right]^{T}$, $r'_{k,0}=\left[0_{1,\tilde{s}_k+1},\,1,\,0_{1,\tilde{s}_k-1}\right]^{T}$, $x'_{k,0}=\left[0_{1,2\tilde{s}_k+1}\right]^T$
\State {{\color{red}$\phi = (\rhat_{k,0}^{'T} G_{k,\tilde{s}_k} \rhat'_{k,0})^{1/2}$}}
\For{$j = 0$ to $\tilde{s}_k-1$}
\State $s_k=j+1$
\State $\alpha_{m+j}=\big({r}'^T_{k,j}G_{k,\tilde{s}_k}r'_{k,j}\big)/\big({p}'^T_{k,j}G_{k,\tilde{s}_k}\mathcal{B}_{k,\tilde{s}_k}p'_{k,j}\big)$ 
\State $q'_{k,j} = \alpha_{m+j}p'_{k,j}$

\State $x'_{k,j+1}=x'_{k,j}+q'_{k,j}$
\State $r'_{k,j+1}=r'_{k,j}-{\mathcal{B}_{k,\tilde{s}_k}}q'_{k,j}$ 

\State $\beta_{m+j}=\big({r}'^T_{k,j+1}G_{k,\tilde{s}_k}r'_{k,j+1}\big)/\big({r}'^T_{k,j}G_{k,\tilde{s}_k}r'_{k,j}\big)$ 

\State $p'_{k,j+1}=r'_{k,j+1}+\beta_{m+j}p'_{k,j}$ 

\State{{\color{red}$\phi = \max \{\phi, (r_{k,j+1}^{'T} G_{k,\tilde{s}_k} r'_{k,j+1} )^{1/2} \}$}}
\State {{\color{red} $\psi_{m+j+1} = \psi_{m+j}/(\psi_{m+j}+\beta_{m+j})$ }}

\If {$m+j+1>1$}
		\State {{\color{red} Update estimates $\lmin$ and $\lmax$.}}
		\State{{\color{red}$c_{m+j+1} = \max\{1,\lmax(\psi_{m+j+1}/\lmin)^{1/2} \}$ }}
\EndIf

\IfThen {{\color{red}$j< \tilde{s}_k-1$ \textbf{and} $\gamma_{j+1} \geq \frac{\veps^*}{ c_{m+j+1} \veps \phi}$ }}{break from inner loop.}

\EndFor

\State {Recover iterates $\{p_{m+s_k},r_{m+s_k},x_{m+s_k}\}$ according to~\eqref{eq:coords}.}
\State {$m=m+s_k$}
\IfThen {$m>1$}{{\color{red}update basis parameters w/ $\lmin$, $\lmax$ by~\eqref{eq:newtbasis} or~\eqref{eq:chebbasis}.}}

\EndFor

\end{algorithmic}
\end{algorithm}

\section{Numerical experiments}
\label{sec:exp}

In this section we present experiments run in MATLAB (version R2017a) to compare the numerical behavior of HSCG, fixed $s$-step CG, adaptive $s$-step CG with a monomial basis, and improved adaptive $s$-step CG with dynamically updated Newton and Chebyshev bases for small SPD matrices from the SuiteSparse collection~\cite{suitesparse}. 
For each matrix, we test the $s$ values $5, 10,$ and $15$ (which are the maximum allowable $s_k$ values $\sigma$ in the adaptive algorithms). We test two different values of $\veps^*$; the first being the relative true residual $2$-norm attainable by HSCG as determined experimentally, and the second being $\veps^* = 10^{-6}$. The adaptive algorithms all use $f=\sigma$ to allow for the largest possible $s_k$ values. 
For all experiments in this section, we use two-sided diagonal preconditioning, where the resulting preconditioned matrix is $ D^{-1/2}AD^{-1/2}$ where $D$ is a diagonal matrix of the largest entries in each row of $A$. The (unpreconditioned) right hand side $b$ is set with entries $1/\sqrt{N}$ and the initial guess is $x_0=0$. We use double precision in all tests ($\veps\approx 2^{-53}$). We discuss the relative performance of the (fixed and adaptive) $s$-step algorithms in terms of the number of outer loop iterations, which can be seen as a proxy for the number of global synchronizations. 
We stress, however, that the experiments here are on very small matrices, and thus the aim is solely to demonstrate the numerical behavior rather than evaluate potential parallel performance improvements.  

\subsection{The benefit of adaptively setting $c_{m+j+1}$}
\label{sec:ck}

\begin{figure}
\centering
\includegraphics[width=2.5in,trim=1.25in 3.2in 1.85in 3.6in, clip]{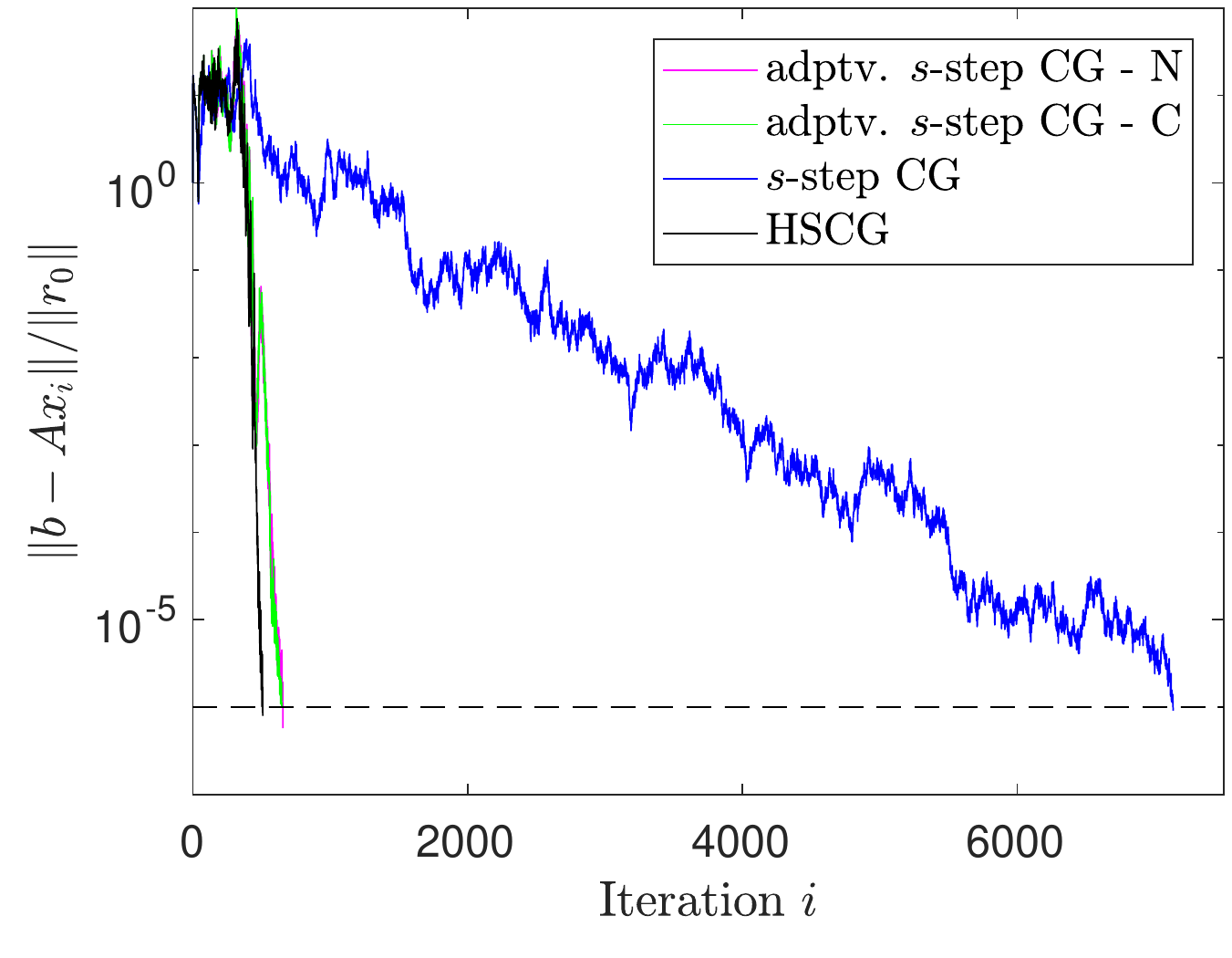} 
\includegraphics[width=2.5in,trim=1.4in 3.2in 1.7in 3.6in, clip]{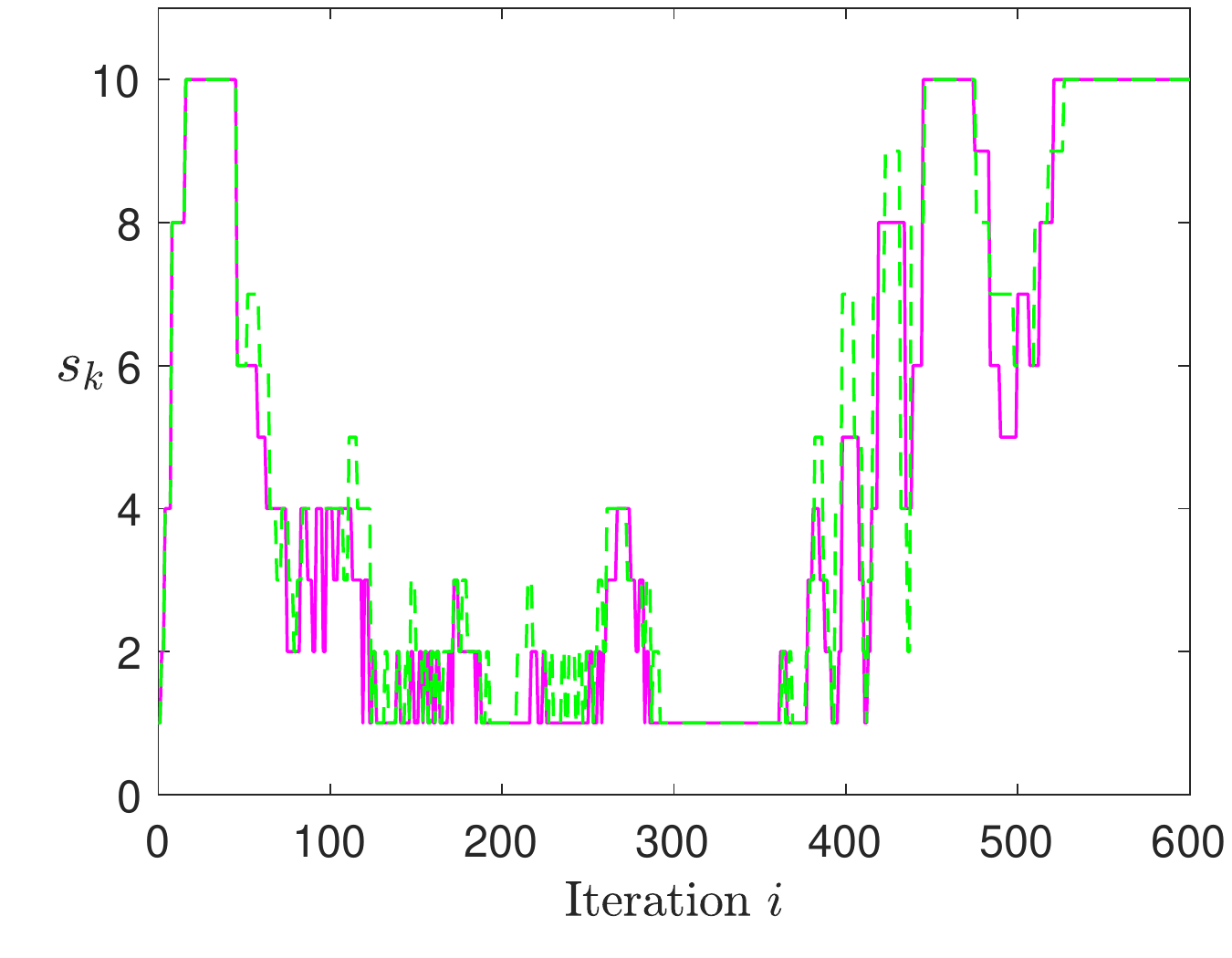} \\
\includegraphics[width=2.5in,trim=1.25in 3.2in 1.85in 3.6in, clip]{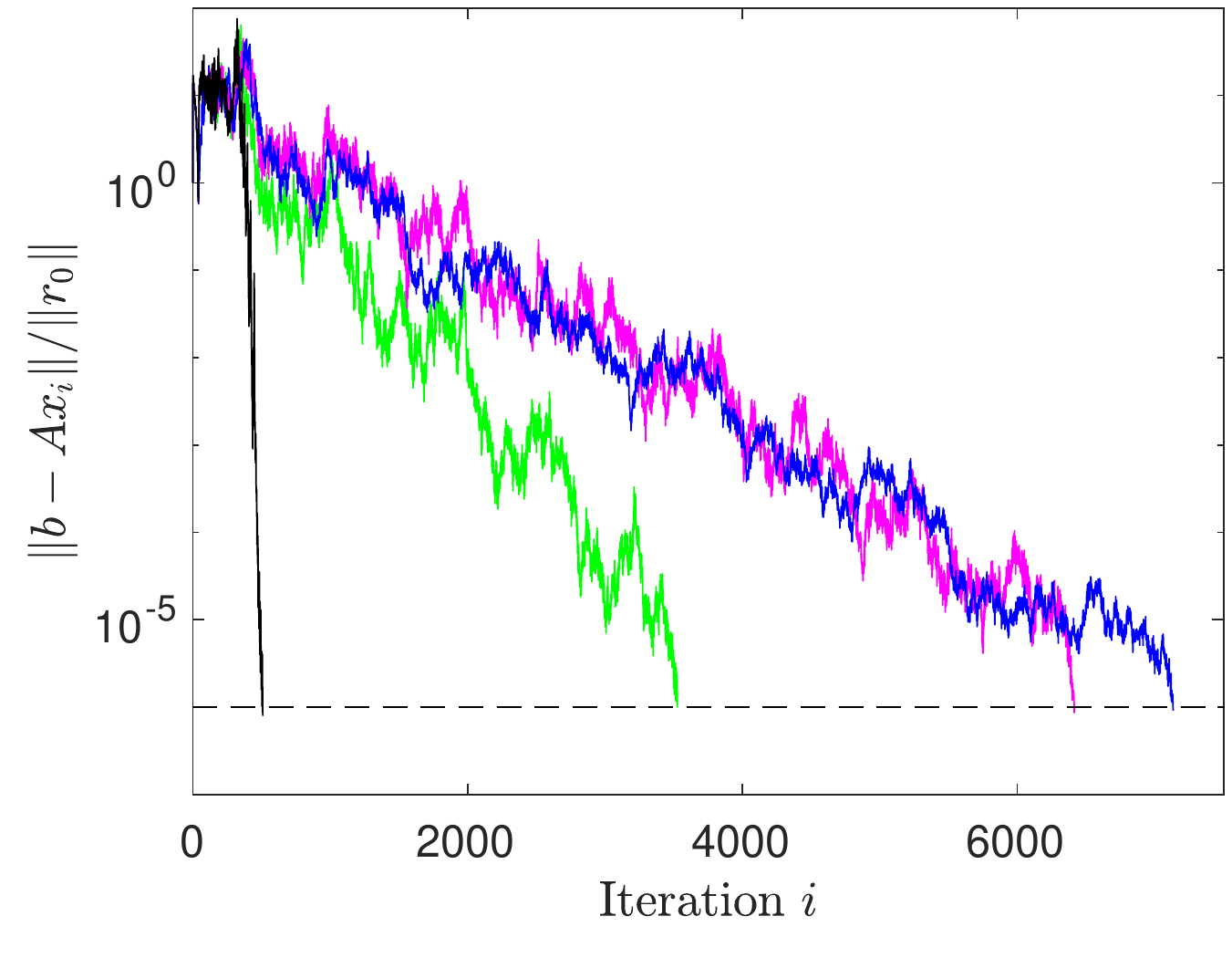}
\includegraphics[width=2.5in,trim=1.4in 3.2in 1.7in 3.6in, clip]{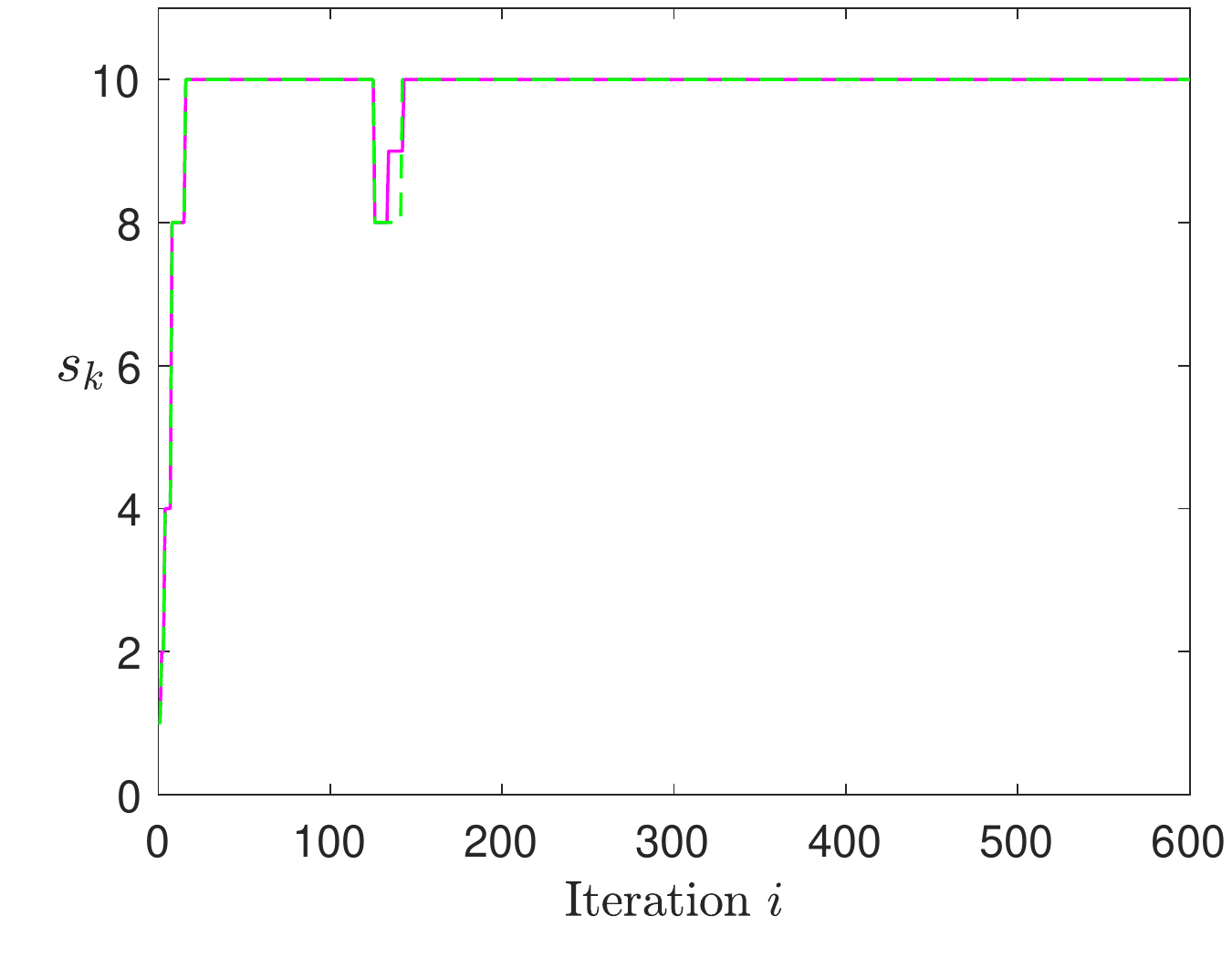} \\
\includegraphics[width=2.5in,trim=1.25in 3.2in 1.85in 3.6in, clip]{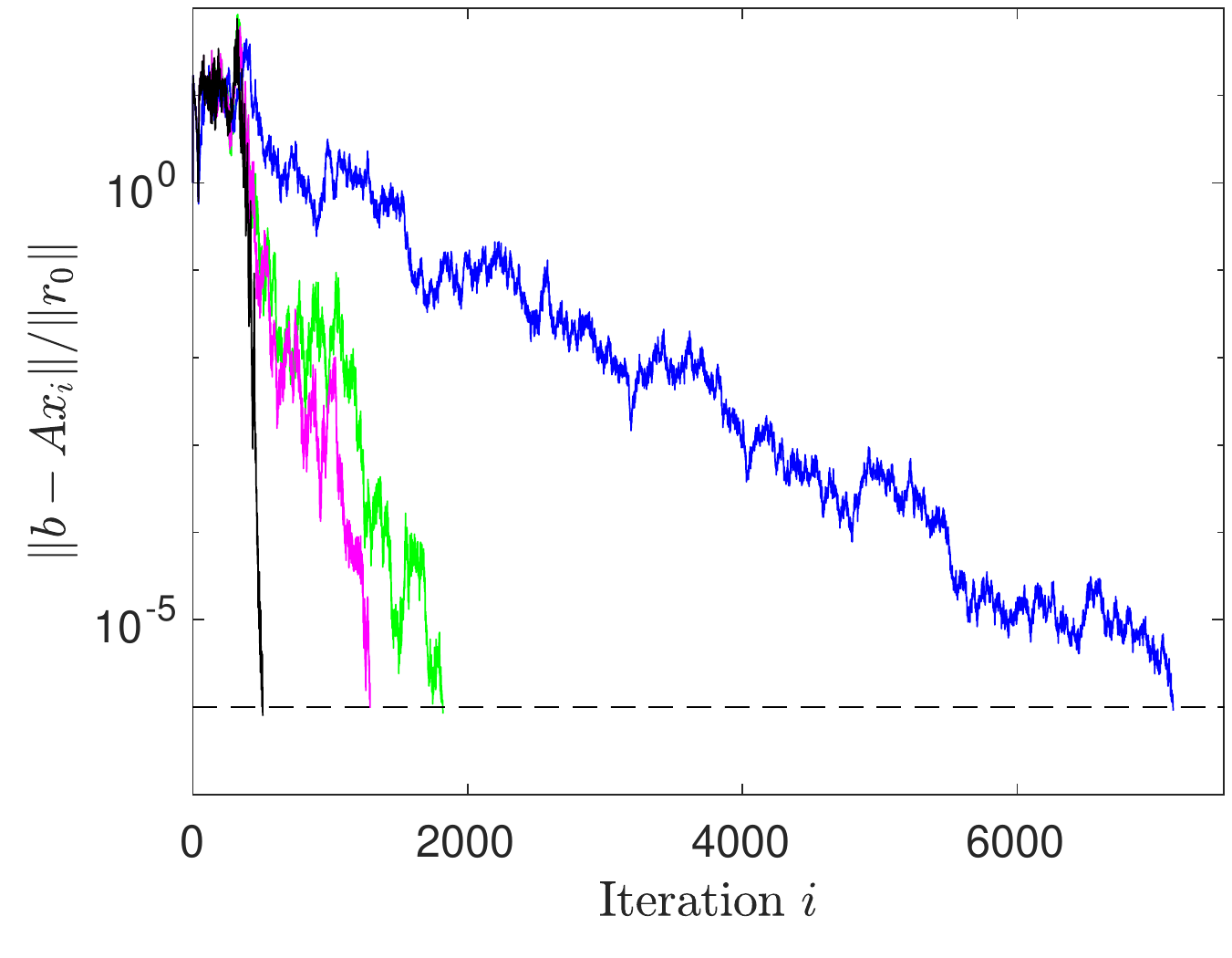} 
\includegraphics[width=2.5in,trim=1.4in 3.2in 1.7in 3.6in, clip]{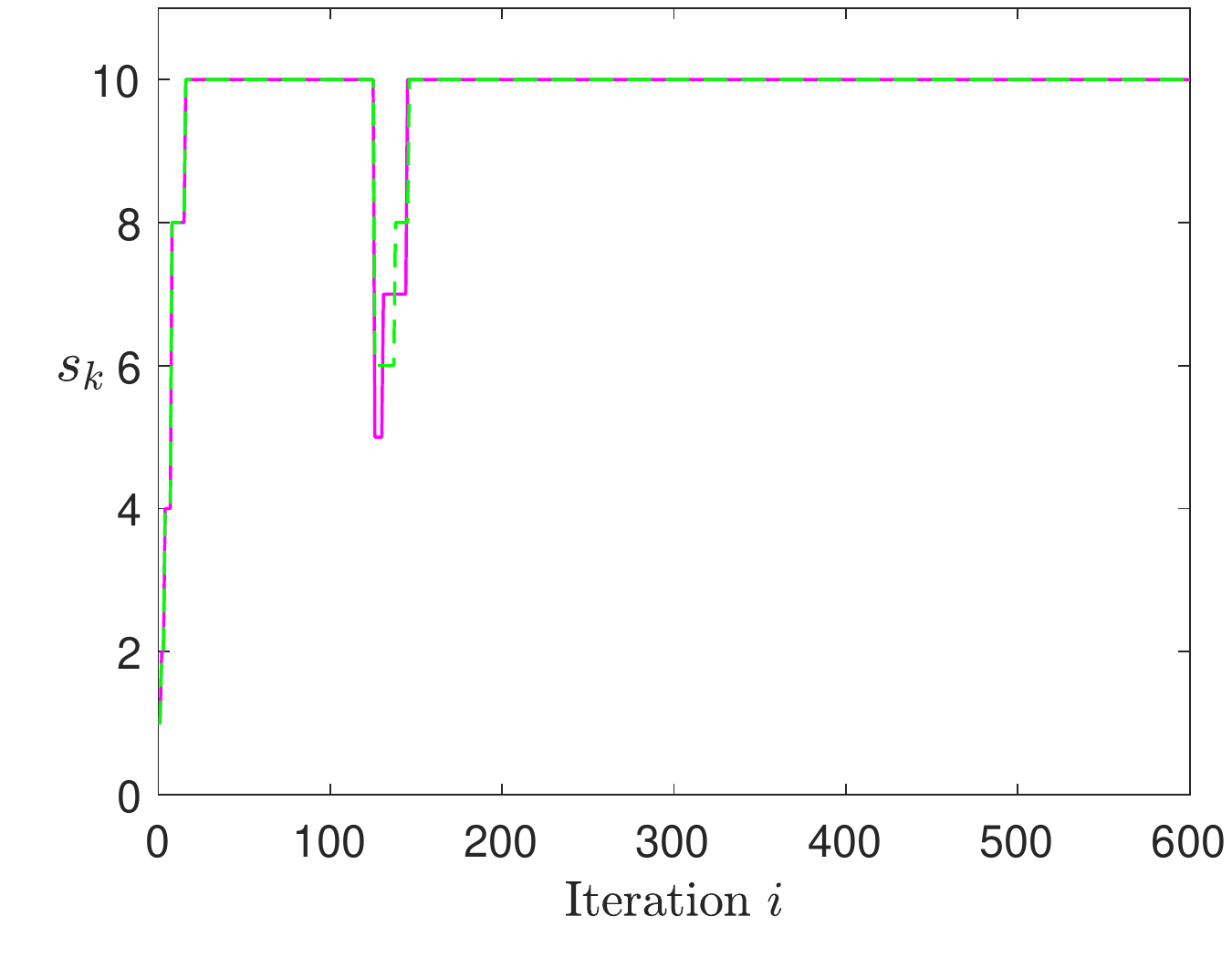}
\caption{\footnotesize{Test results for matrix \texttt{nos1} with $s,\sigma=10$ and $\veps^*=10^{-6}$. Plots on the left show the convergence trajectories for HSCG, fixed $s$-step CG, and improved adaptive $s$-step CG with Newton and Chebyshev bases. Plots on the right show the 
values of $s_k$ in the improved adaptive algorithm using the Newton basis (magenta) and the Chebyshev basis (green). The top row uses $c_{m+j+1}=\lmax/\lmin$, the middle row uses $c_{m+j+1}=1$, and the bottom row uses the new adaptive approach, setting $c_{m+j+1}=\tilde{\xi}_{m+j+1}$; see~\eqref{eq:xiupd}.}}
\label{fig.nos1test}
\end{figure}

We first demonstrate the benefits of our approach for dynamically setting the constant $c_{m+j+1}$. 
In Figure~\ref{fig.nos1test}, we plot the convergence (in terms of relative residual $2$-norm) of HSCG, fixed $s$-step CG, and the improved adaptive $s$-step CG algorithm (Algorithm~\ref{alg:vscg2}) for the matrix \texttt{nos1} from SuiteSparse~\cite{suitesparse} using various values of $c_{m+j+1}$, with the specified tolerance $\varepsilon^* = 10^{-6}$ and $s=10$ ($\sigma=10$ in the adaptive algorithm). Properties of the \texttt{nos1} matrix can be found in Table~\ref{matrices}. For this problem, HSCG takes 510 iterations to converge. Fixed $s$-step CG suffers from significantly delayed convergence, requiring a total of 7134 iterations (corresponding to 714 outer loop iterations) to converge. We note that the old adaptive $s$-step CG (Algorithm~\ref{alg:vscg}) in this case chooses $s_k=\sigma$ in each outer loop iteration and thus also requires $714$ outer iterations to converge. Thus for these algorithms we would expect performance to be worse than HSCG. 

For the improved adaptive $s$-step CG algorithm, we show convergence for 3 different choices of $c_{m+j+1}$. The plots on the left show convergence of the relative true residual and the plots on the right show the $s_k$ values used throughout the iterations. 
The top row of plots uses $c_{m+j+1}=\lmax/\lmin$ where $\lmax$ and $\lmin$ are updated in every iteration, the middle row uses a constant $c_{m+j+1} = 1$, and the bottom row uses the new adaptive approach which automatically sets $c_{m+j+1}=\tilde{\xi}_{m+j+1}$; see~\eqref{eq:xiupd}. 

In the top row of plots using $c_{m+j+1}=\lmax/\lmin$, the convergence of the improved adaptive $s$-step CG algorithm is relatively close to that of the HSCG algorithm for both Newton and Chebyshev bases. However, this comes at the cost of the use of small $s_k$ values, with $s_k=1$ in many outer loop iterations (see the top right plot). Here the improved adaptive $s$-step CG algorithm requires 291 and 255 outer iterations with the Newton and Chebyshev bases, respectively. 

It is clear from the middle plots that blindly using $c_{m+j+1}=1$ as previously suggested in~\cite{carson18} is not always sufficient in practice (note the stretched $x$-axis in the left-hand plot for this experiment). Using $c_{m+j+1}=1$ results in larger $s_k$ values (with the maximum $s_k=10$ occurring in most outer loop iterations), but this results in a significant convergence delay due to the use of ill-conditioned bases. This results in a \emph{greater} total number of outer loop iterations required for convergence, with the implementations with Newton and Chebyshev bases requiring 646 and 357 outer loop iterations, respectively. Note that for the Newton basis case, this is more ``global synchronizations'' than in HSCG!

\begin{figure}
\centering
\includegraphics[width=3in,trim=1.25in 3.2in 1.25in 3.6in, clip]{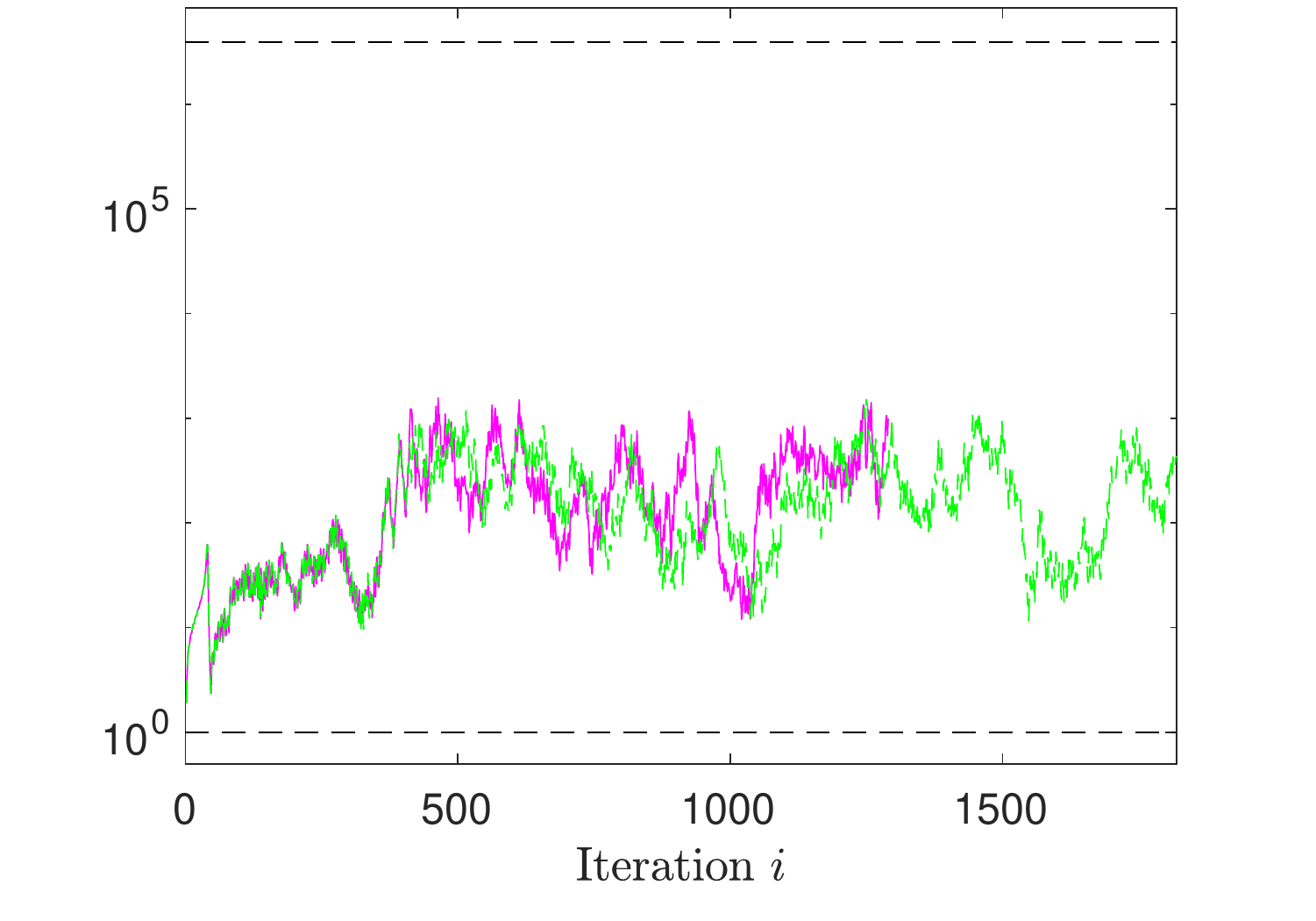}
\caption{\footnotesize{The value of $c_{m+j+1}=\tilde{\xi}_{m+j+1}$ for improved adaptive $s$-step CG using the Newton basis (magenta) and the Chebyshev basis (green). The dashed horizonal lines show the bounds for $\xi_{m+j+1}$, $1$ and $\kappa(A)$.}}
\label{fig.nos1testc}
\end{figure}

Finally, in the bottom row, we see the advantage of our approach for automatically setting $c_{m+j+1}$ based on information gained from the iterations. Here the total number of iterations required for convergence is greater than when we use $c_{m+j+1}=\lmax/\lmin$, but still much smaller than the naive choice $c_{m+j+1}=1$. But notice that the $s_k$ values used are much closer to the $c_{m+j+1}=1$ case, with the maximum value $s_k=10$ occurring in the majority of iterations. This results in fewer overall outer loop iterations; here the implementations with the Newton and Chebyshev bases required only 134 and 187 outer loop iterations, respectively. For reference, in Figure~\ref{fig.nos1testc}, we show the value of $c_{m+j+1}=\tilde{\xi}_{m+j+1}$ computed in each iteration. The dashed horizontal lines mark 1 (the minimum value of $\xi_{m+j+1}$) and $\kappa(A)$ (the maximum value of $\xi_{m+j+1}$). We see that the computed $\tilde{\xi}_{m+j+1}$ is actually somewhere between $1$ and $\kappa(A)$; by using this tighter, computable bound, we effectively balance the goals of minimizing the number of iterations required for convergence and maximizing the $s_k$ values used.

\subsection{Test problems from SuiteSparse}
\label{sec:ss}

We now present a few select test problems which demonstrate the behavior of the improved adaptive $s$-step CG algorithm. For all problems, the experimental setup is as described at the beginning of Section~\ref{sec:exp}. The matrices used in the experiments are shown in Table~\ref{matrices}, where the norm and condition number are those of the diagonally preconditioned system matrices.  For space purposes, we only include plots for two of the test matrices, \texttt{494bus} and \texttt{bcsstk09}. However, for each test matrix, we include a table which shows the number of iterations for HSCG and the number of outer loop iterations (global synchronizations) for fixed $s$-step CG, adaptive $s$-step CG, and improved adaptive $s$-step CG with both Newton and Chebyshev bases, for the $s$ (or $\sigma$) values 5, 10, and 15, and for two different choices of $\veps^*$. For the $s$-step variants, the first number gives the number of outer loop iterations and the number inside the parentheses gives the total number of (inner loop) iterations. A dash alone in the table indicates that the algorithm diverged. A dash along with a number in square brackets indicates that the algorithm did eventually converge (i.e., the residual norm stagnated), but to a relative residual norm less than $\veps^*$; the number in square brackets gives the final attained relative residual $2$-norm. 
In the plots, markers are used to denote the outer iterations in the (fixed and adaptive) $s$-step algorithms, which are the points at which global synchronization would occur.  

\begin{table}[hbt!]
\centering
\footnotesize
\caption{\footnotesize{Test matrix properties}}
\label{matrices}
\begin{tabular}{|c|c|c|c|c|}
\hline
Matrix & $N$ & nnz & $\Vert A \Vert$ & $\kappa(A)$ \\ \hline
494bus & 494 & 1666 & 2.00 & $7.90\cdot 10^4$ \\ \hline
bcsstk09 & 1083 & 18437 & 1.98 & $1.04\cdot 10^4$ \\ \hline
gr\_30\_30 & 900 & 7744 & 1.49 & $1.95\cdot 10^2$ \\ \hline
nos6 & 675 & 3255 & 2.00 & $3.49\cdot 10^6$ \\ \hline
mhdb416 & 416 & 2312 & 1.52 & $3.02\cdot 10^3$ \\ \hline
nos1 & 237 & 1017 & 2.00 & $3.94\cdot 10^6$ \\ \hline
\end{tabular}
\end{table}

We show results for the matrix \texttt{494bus} in Figure~\ref{fig.494busn_0} and the corresponding Table~\ref{tab.494busn0}. This problem represents the ideal case for the improved adaptive $s$-step algorithm and nicely highlights the benefits of the improved approach. When $\veps^*=2.2\cdot 10^{-10}$, the improved algorithm outperforms the fixed $s$-step algorithm and the old adaptive $s$-step algorithm in all cases, even for small $\sigma$. When $\veps^*=10^{-6}$, the improved approach takes about the same number of inner and outer loop iterations as $s$-step CG and old adaptive $s$-step CG for $s,\sigma=5$, but clearly outperforms both algorithms for higher $s, \sigma$ values. 

Using the improved bases generated using $\lmin$ and $\lmax$ has a clear benefit; in all cases, the total number of iterations required by the improved adaptive approach is about the same as for HSCG. In contrast, for $\veps^*=2.2\cdot 10^{-10}$, the fixed $s$-step and old adaptive approaches require more than 200 more iterations to converge compared to HSCG even for $s,\sigma=5$. For larger $s,\sigma$ values, fixed $s$-step CG no longer converges to the specified tolerance. Old adaptive $s$-step CG eventually converges for $\sigma=10$, but requires more than twice the total number of iterations as HSCG, limiting the potential benefit of any $s$-step approach. For $s,\sigma=15$, both fixed $s$-step CG and old adaptive $s$-step CG diverge. In the case $\veps^*=10^{-6}$, the number of iterations required for convergence in old adaptive $s$-step CG grows drastically with $\sigma$, more than doubling from $\sigma=5$ to $\sigma=10$ and more than quadrupling from $\sigma=10$ to $\sigma=15$! For $\sigma=15$, this results in a greater number of outer loops iterations (i.e., more global synchronizations) than HSCG! In contrast, the total number of iterations required for the improved approach stays constant with $\sigma$ also in the $\veps^*=10^{-6}$ case, and the number of outer loop iterations required also continues to decrease with increasing $\sigma$. For $\sigma=15$, the improved approach with both Newton and Chebyshev bases exhibits a decrease in outer loop iterations of more than $12\times$ versus HSCG.

We next show results for the text problem with matrix \texttt{bcsstk09} (Figure~\ref{fig.bcsstk09n_0} and Table~\ref{tab.bcsstk09n0}). For both tested $\veps^*$ values, when $s,\sigma=5$ the improved adaptive approach requires more outer iterations to converge than both fixed $s$-step CG and old adaptive $s$-step CG. However, again the benefit of the more well-condition bases is clear for larger $s,\sigma$ values; as $s,\sigma$ increases, the total number of outer iterations actually increases rather than decreases for both fixed $s$-step CG and old adaptive $s$-step CG. Neither of those approaches converges to the prescribed level $\veps^*=1.9\cdot 10^{-12}$ when $s,\sigma=15$. The old adaptive approach eventually converges to the prescribed level when $\veps^*=10^{-6}$ and $\sigma=15$, although it requires almost $10\times$ the number of total iterations as HSCG; in fact, looking at how the total number of outer loop iteration required grows as $\sigma$ is increased, we see around a $2\times$ increase in the data movement cost for each increase in $\sigma$. In contrast, the number of total iterations required for the improved approach is constant with $\sigma$ using both Newton and Chebyshev bases, and thus the number of outer loop iterations  decreases with increasing $\sigma$ as desired.

For \texttt{gr\_30\_30}, \texttt{nos6}, and \texttt{mhdb416} (Tables~\ref{tab.gr3030n0}, ~\ref{tab.nos6n0}, and~\ref{tab.mhdb416n0}, respectively), the story is similar; fixed $s$-step CG and old adaptive $s$-step CG fail to converge to the requested accuracy in many cases, especially for larger $s,\sigma$ values. In cases where these algorithms do converge for higher $s,\sigma$ values, the number of total iterations required (and also the number of outer loop iterations required) generally grow with increasing $s,\sigma$, sometimes drastically. For example, for the tests with \texttt{mhdb416} for $\veps^*=10^{-6}$ and $\sigma=15$, the number of total iterations required for convergence grows by almost $40\times$ versus HSCG, and thus requires more than $3\times$ the number of global synchronizations versus HSCG. In all cases the improved adaptive approach converges to the prescribed level and in most cases the number of iterations required for convergence remains relatively constant, resulting in a decrease in the number of outer loop iterations with increasing $\sigma$. We note that in some cases, the total number of outer loop iterations stagnates with increasing $\sigma$, or at least the rate of decrease is decreasing; see, for example, the experiments with \texttt{gr\_30\_30} with $\sigma =10,15$ for both $\veps^*$ values. Here, there is no decrease in total number of outer loop iterations going from $\sigma=10$ to $\sigma=15$ for either $\veps^*$. This behavior can occur due to some fundamental limit on how large $s$ can be before the bases become too ill-conditioned. Regardless, this demonstrates another benefit of the adaptive approach; even if the user inputs (or an autotuner chooses) a $\sigma$ value that is too large, the adaptive approach can still be successful.


The experiments with the matrix \texttt{nos1} (Table~\ref{tab.nos1n0}, see also Section~\ref{sec:ck}) are a notable exception. This is a particularly difficult problem even for HSCG, which requires more than $2N$ iterations to converge. Even for $\sigma=5$, the improved adaptive algorithm requires more than $2\times$ the number of iterations to converge than HSCG. We are not sure why this behavior occurs. Regardless, the benefits of the improved approach are still clearly observed. Neither fixed $s$-step 
CG nor old adaptive $s$-step CG converges to the prescribed accuracy when $s,\sigma=15$ is used, and even when convergence does occur, the number of outer loop iterations required grows quickly with $s,\sigma$.  In all cases, the improved adaptive approach still converges to the prescribed level and still reduces the number of outer iterations required by at least a factor of $2$ versus HSCG. We also see that the number of outer loop iterations required in the improved approach still reliably decreases with increasing $\sigma$.

\afterpage{
\begin{figure}[htt!]
\centering
\vspace{-2mm}
\includegraphics[width=2.5in,trim=1.25in 3.2in 1.85in 3.15in, clip]
	{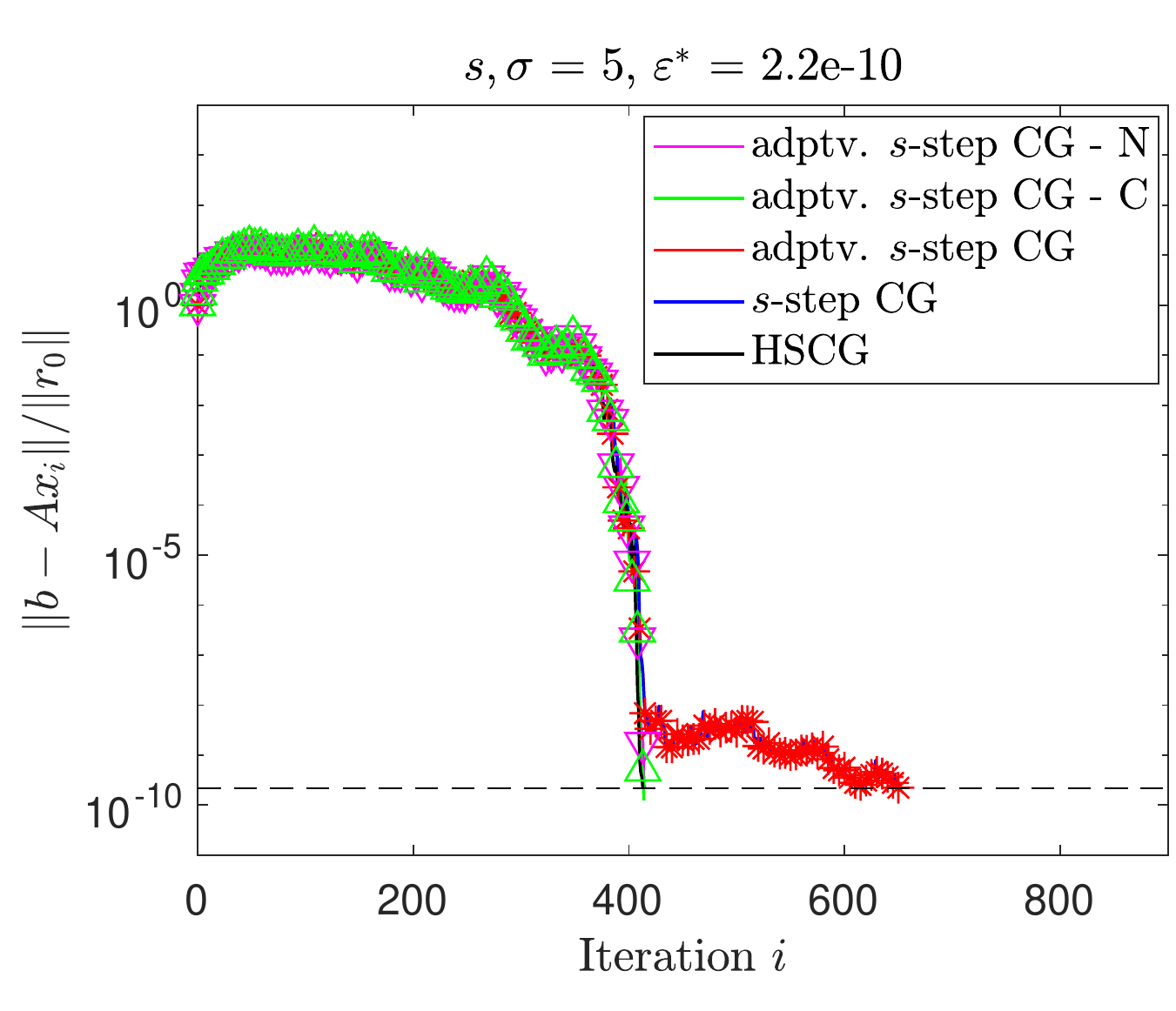}
\includegraphics[width=2.5in,trim=1.3in 3.2in 1.8in 3.15in, clip]
	{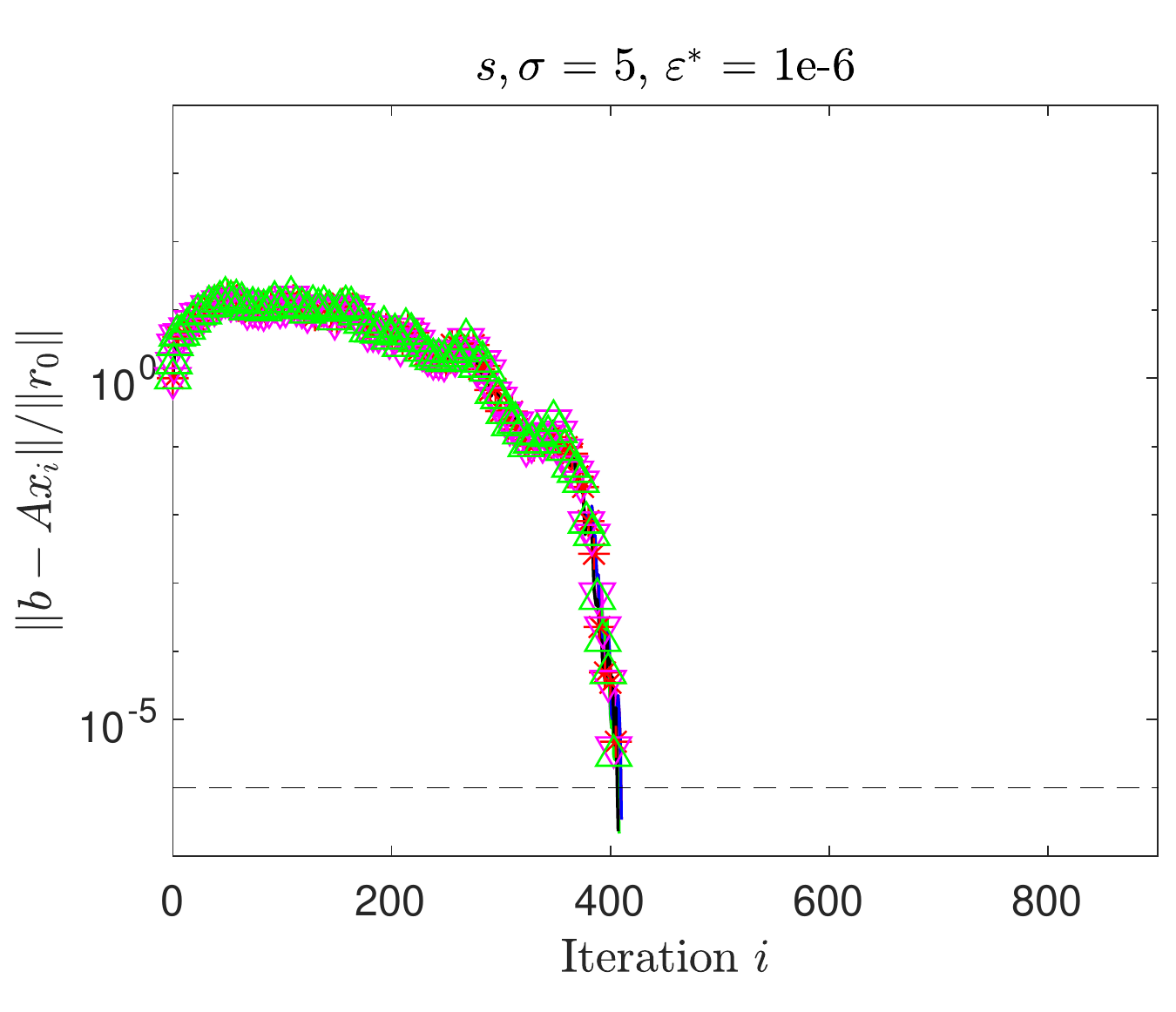}  \\ 
\includegraphics[width=2.5in,trim=1.25in 3.2in 1.85in 3.15in, clip]
	{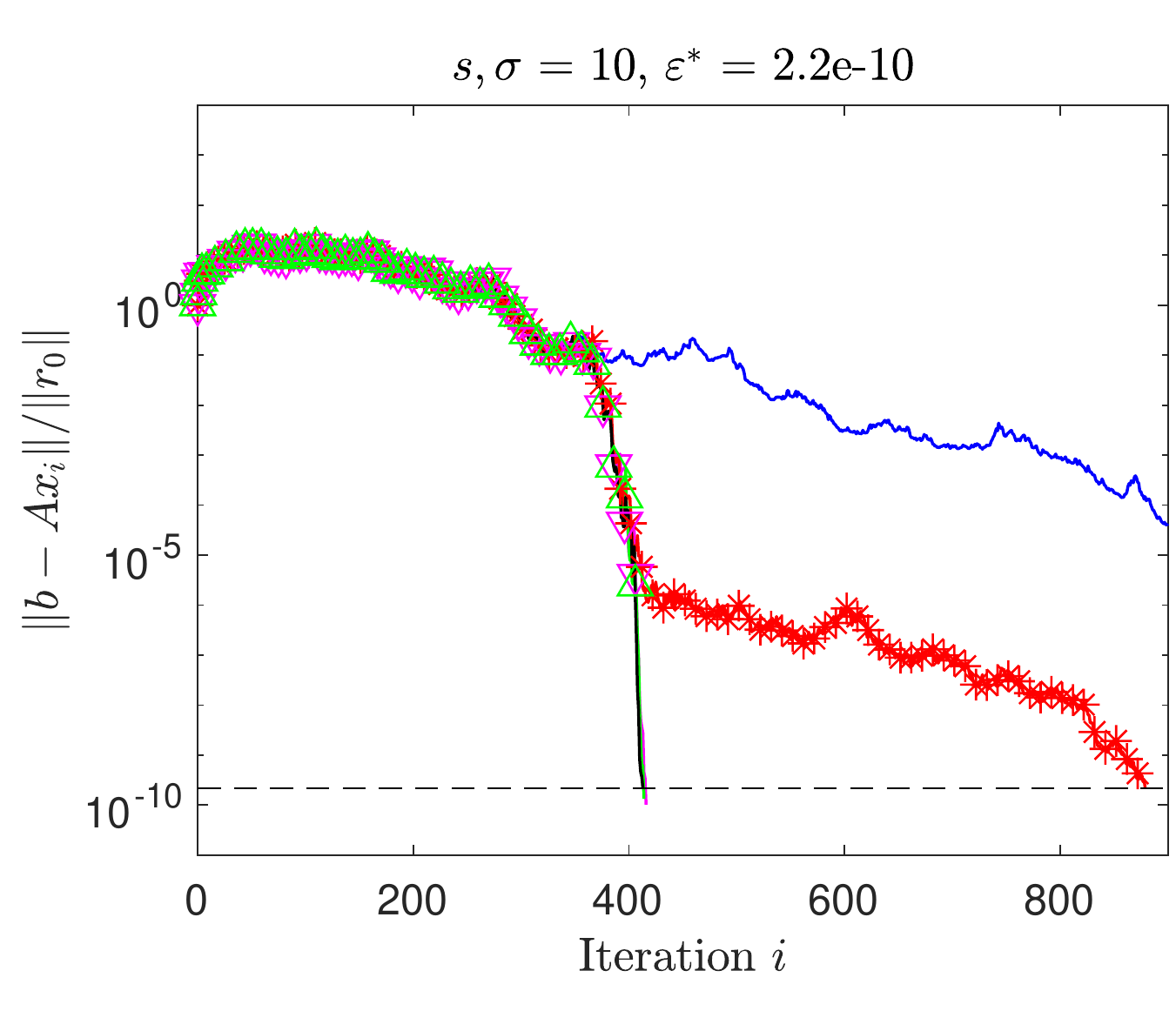}
\includegraphics[width=2.5in,trim=1.3in 3.2in 1.8in 3.15in, clip]
	{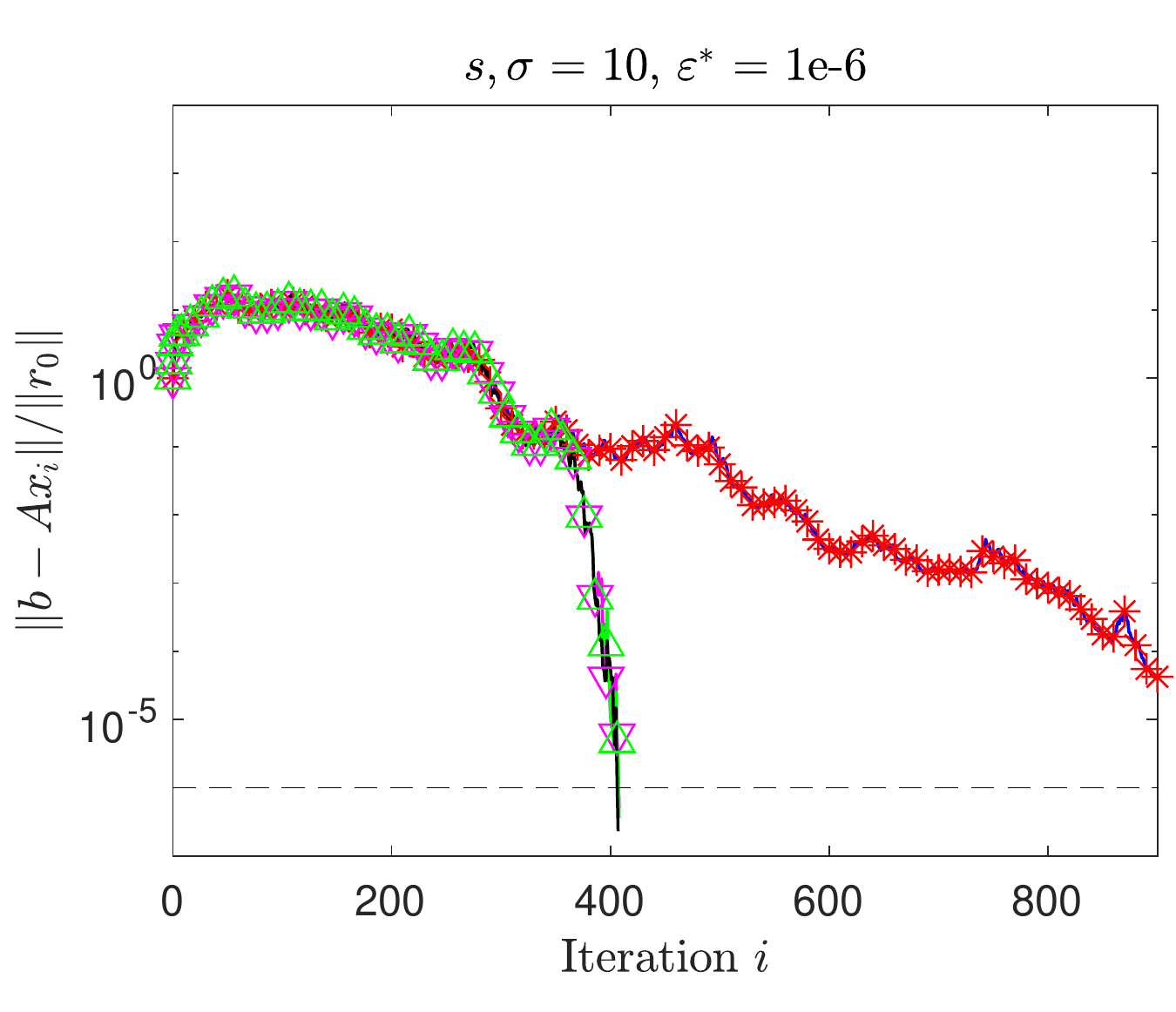} \\ 
	\includegraphics[width=2.5in,trim=1.25in 3.2in 1.85in 3.15in, clip]
	{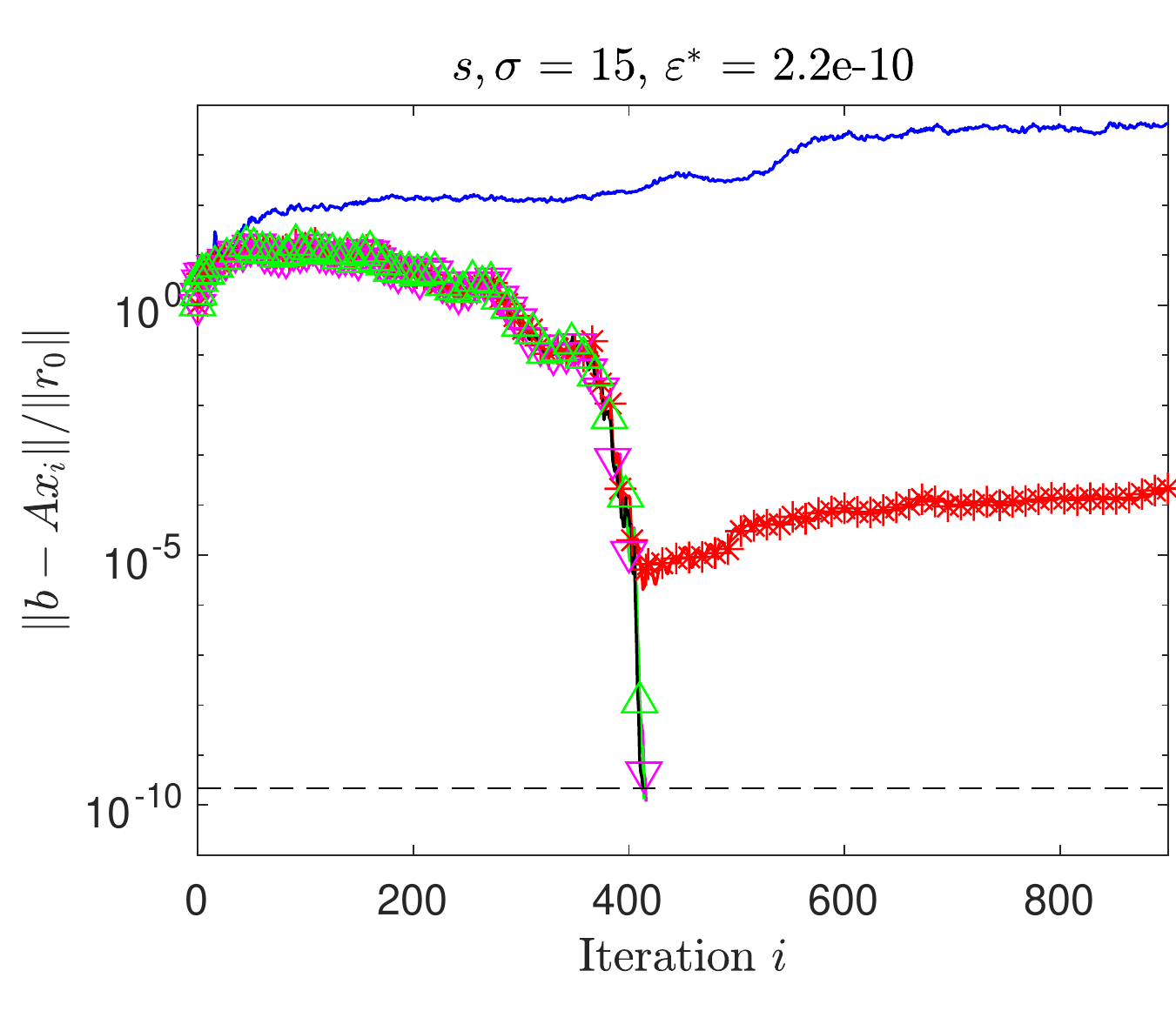}
\includegraphics[width=2.5in,trim=1.3in 3.2in 1.8in 3.15in, clip]
	{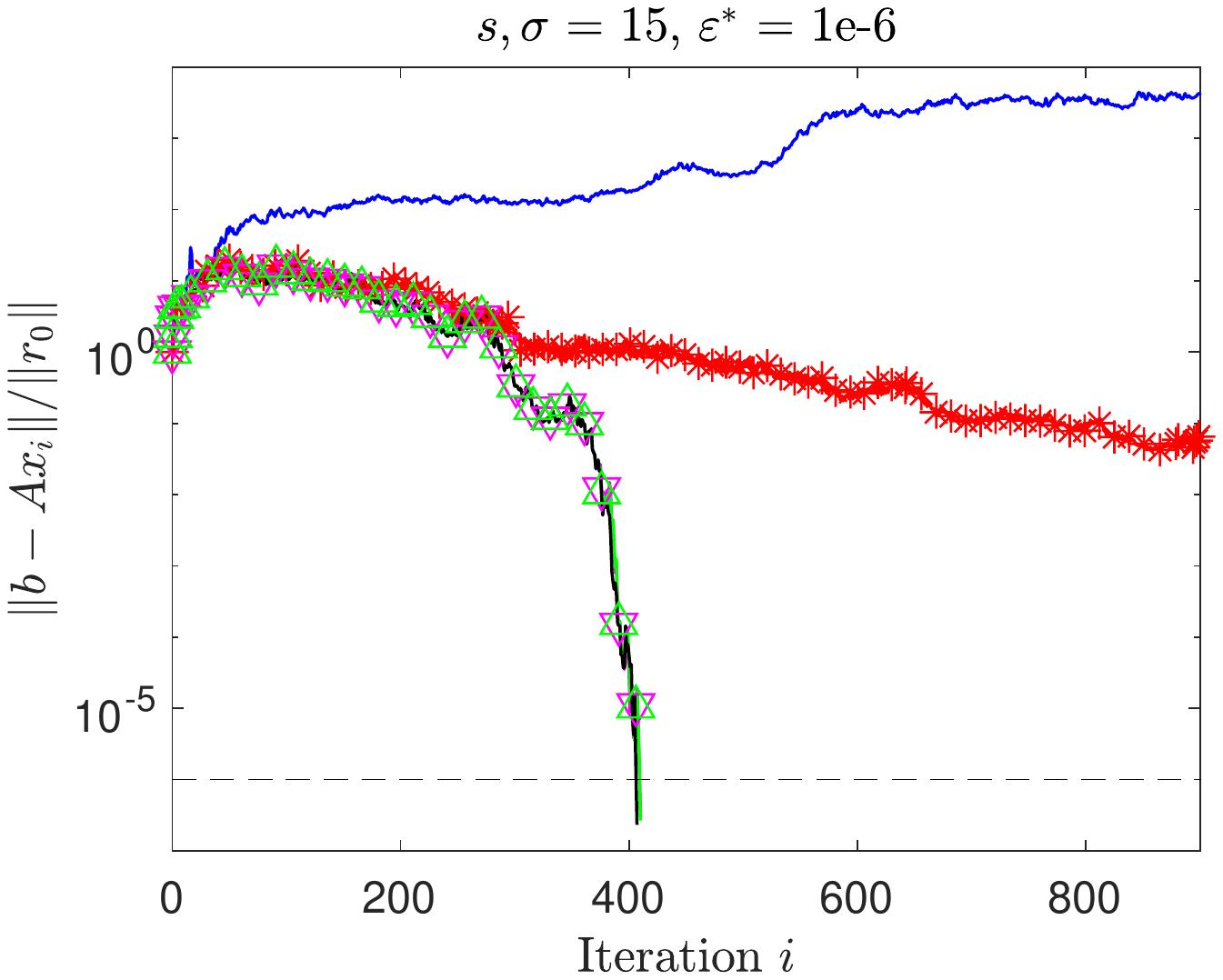} 
\caption{\footnotesize{Convergence of the relative true residual 2-norm for the matrix \texttt{494bus}. 
}}
\label{fig.494busn_0}
\end{figure}

\begin{table}[h!]
\centering
\resizebox{\textwidth}{!}{%
\begin{tabular}{ll|c|c|c|c|c|}
\cline{3-7}
 &  & fixed s-step CG & old adptv s-step CG & \begin{tabular}[c]{@{}c@{}}impr. adptv s-step\\ CG w/Newton basis\end{tabular} &\begin{tabular}[c]{@{}c@{}}impr. adptv s-step\\ CG w/Chebyshev basis\end{tabular} & HSCG \\ \hline
\multicolumn{1}{|l|}{\multirow{3}{*}{$\varepsilon^*=$ 2.2e-10}} & $s,\sigma=5$ & 131 (652) & 131 (652) & 86 (415) & 86 (414) & \multirow{3}{*}{413} \\ \cline{2-6}
\multicolumn{1}{|l|}{} & $s,\sigma=10$ & - [2e-08] & 109 (879) & 58 (416) & 53 (414) &  \\ \cline{2-6}
\multicolumn{1}{|l|}{} & $s,\sigma=15$ & - & - & 57 (416) & 51 (414) &  \\ \hline
\multicolumn{1}{|l|}{\multirow{3}{*}{$\varepsilon^*=$ 1e-6}} & $s,\sigma=5$ & 82 (410) & 82 (410) & 84 (408) & 84 (408) & \multirow{3}{*}{407} \\ \cline{2-6}
\multicolumn{1}{|l|}{} & $s,\sigma=10$ & 115 (1147) & 115 (1147) & 45 (408) & 45 (408) &  \\ \cline{2-6}
\multicolumn{1}{|l|}{} & $s,\sigma=15$ & - & 442 (6146) & 32 (410) & 32 (410) &  \\ \hline
\end{tabular}
}

\caption{\footnotesize{Results for experiments with the matrix \texttt{494bus}. 
}}
\label{tab.494busn0}
\end{table}
}


\afterpage{
\begin{figure}[htt!]
\centering
\includegraphics[width=2.5in,trim=1.25in 3.2in 1.85in 3.15in, clip]
	{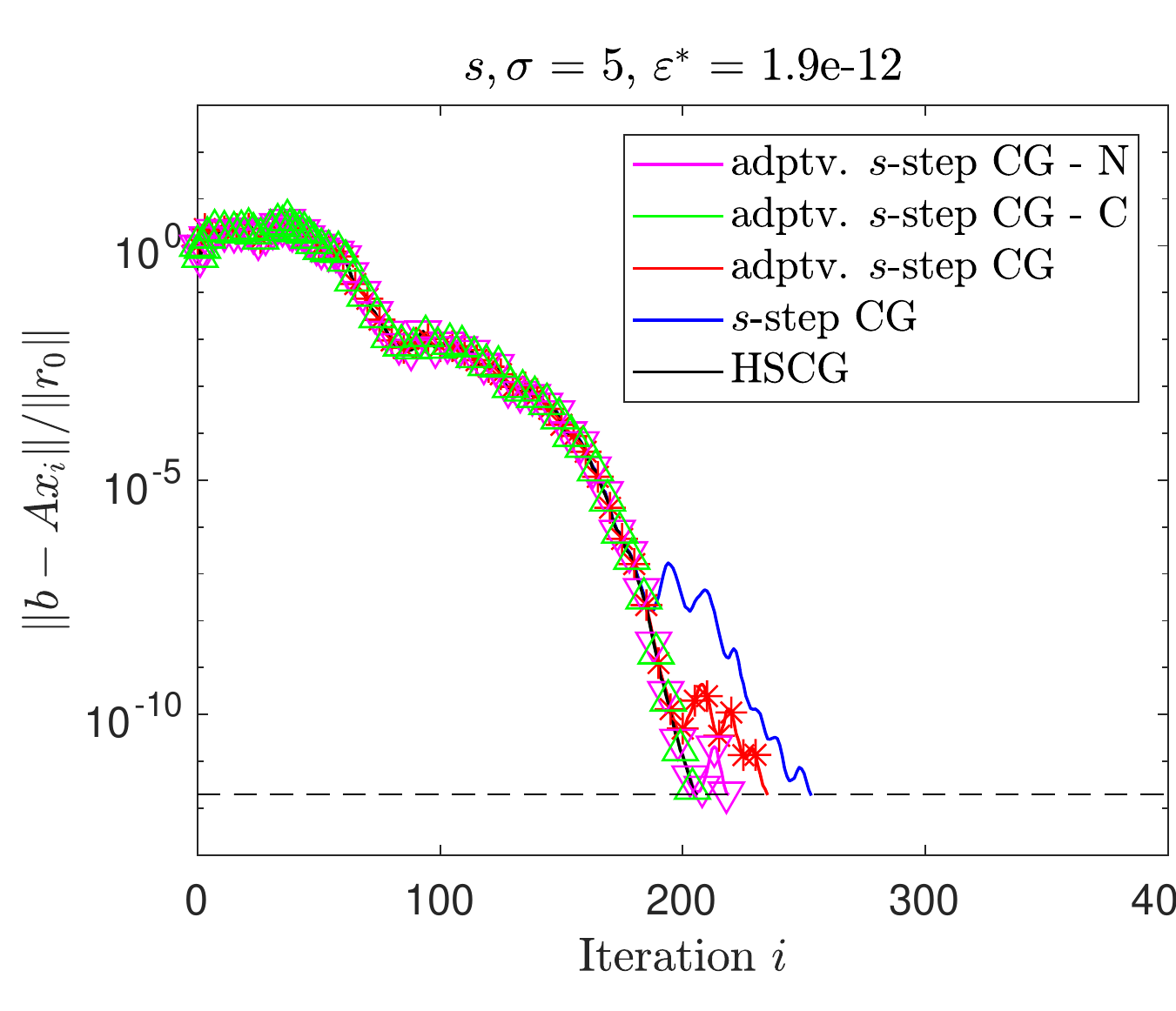}
\includegraphics[width=2.5in,trim=1.25in 3.2in 1.85in 3.15in, clip]
	{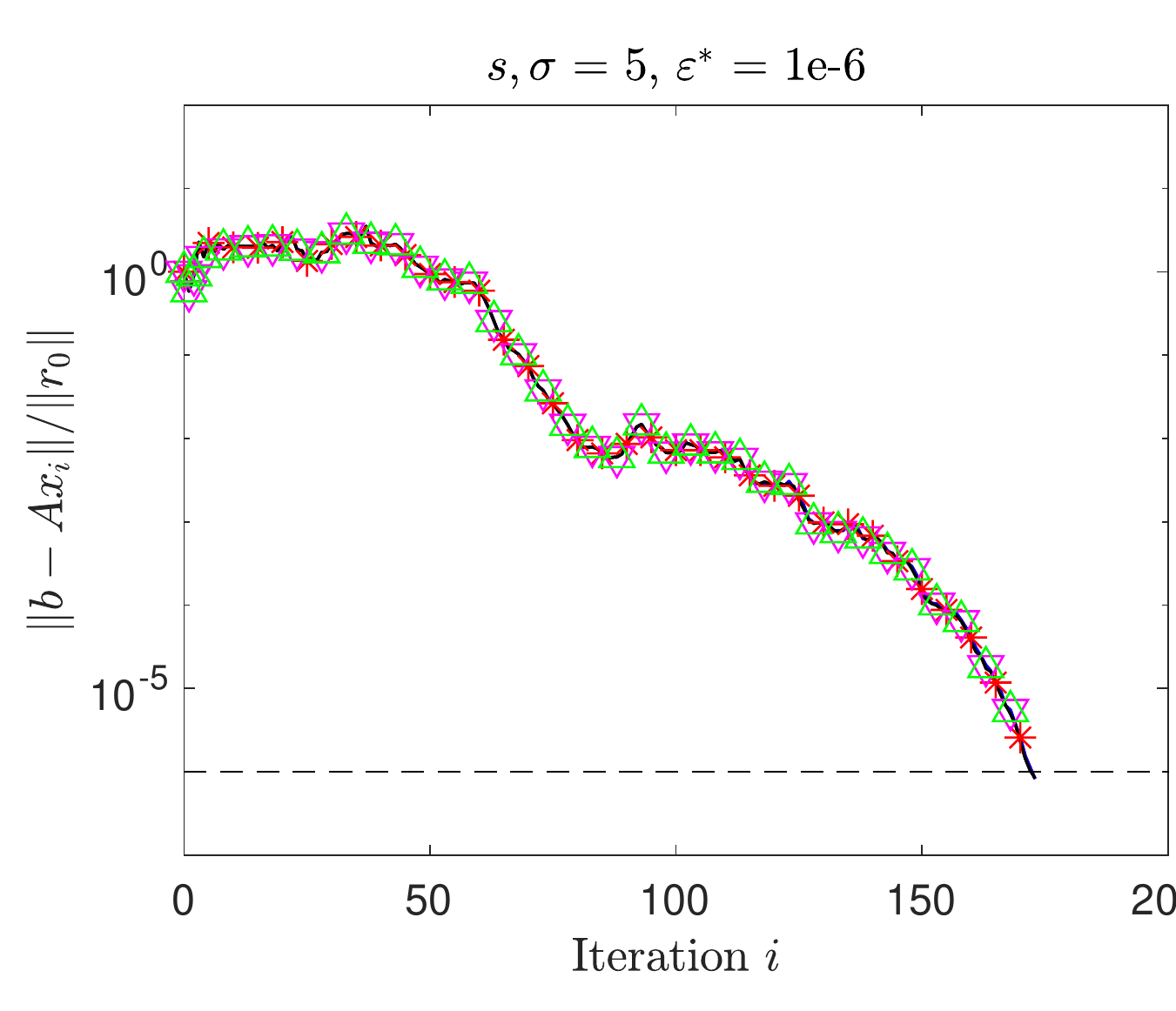} \\
\includegraphics[width=2.5in,trim=1.25in 3.2in 1.85in 3.15in, clip]
	{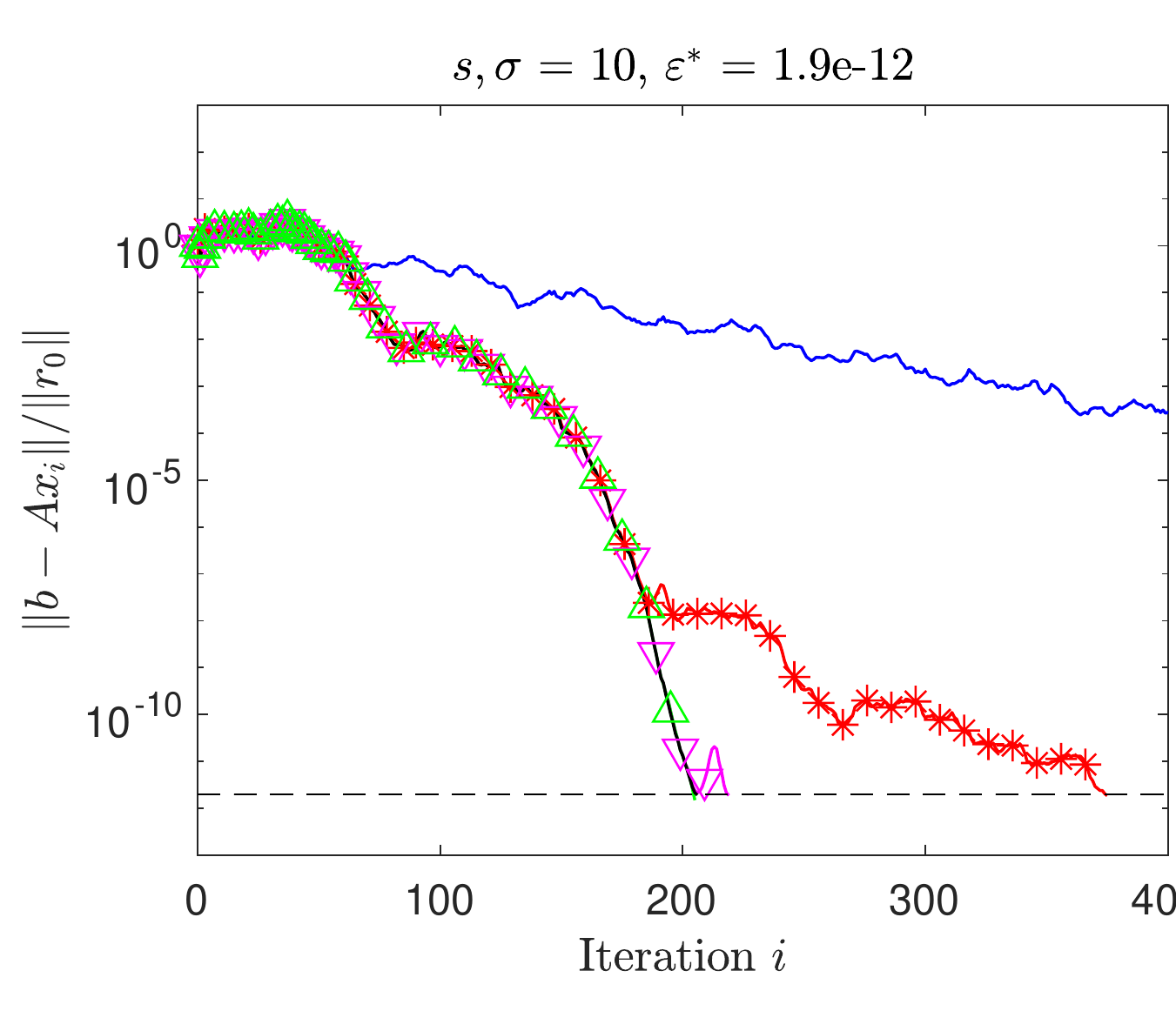}
\includegraphics[width=2.5in,trim=1.25in 3.2in 1.85in 3.15in, clip]
	{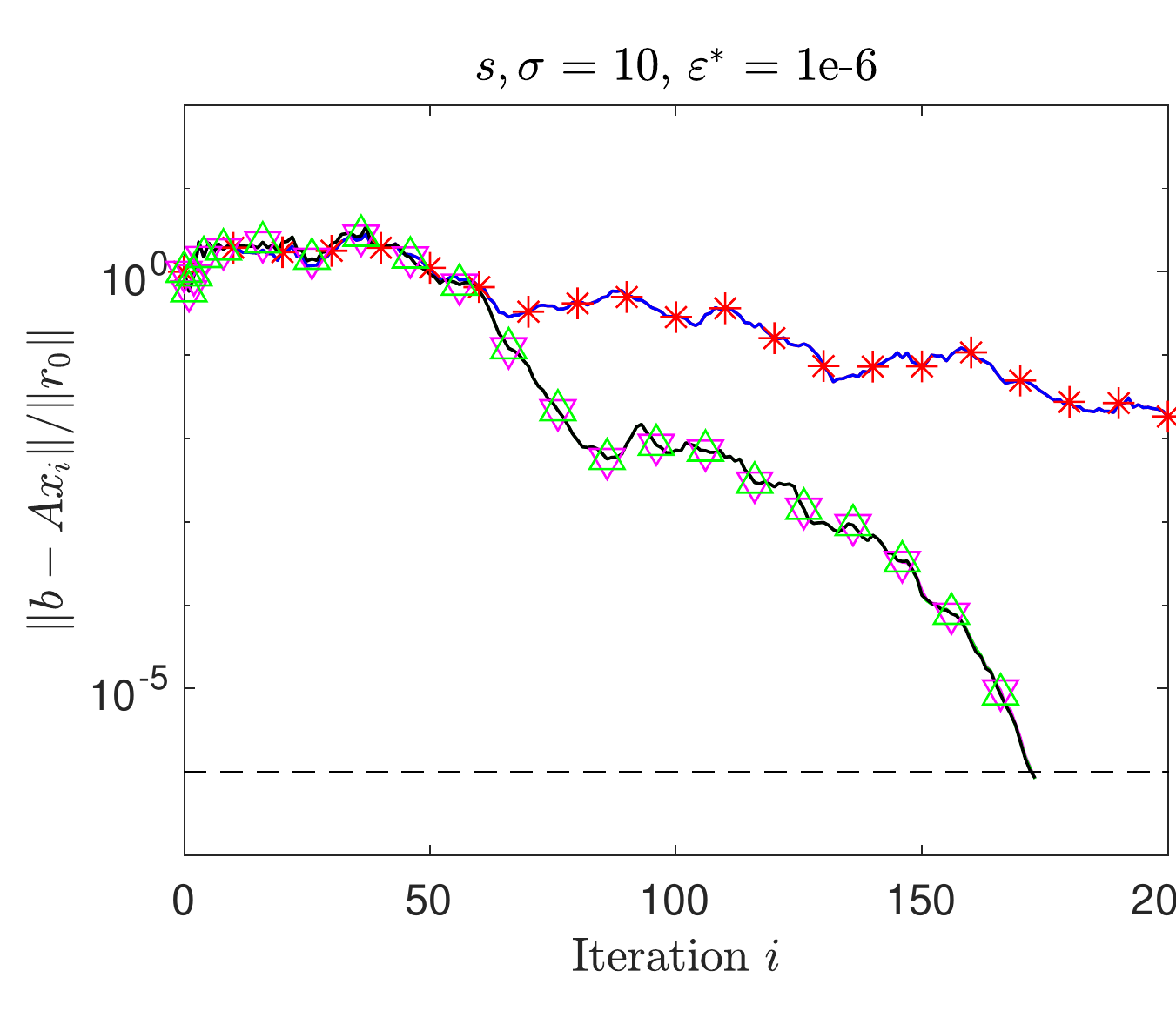} \\
	\includegraphics[width=2.5in,trim=1.25in 3.2in 1.85in 3.15in, clip]
	{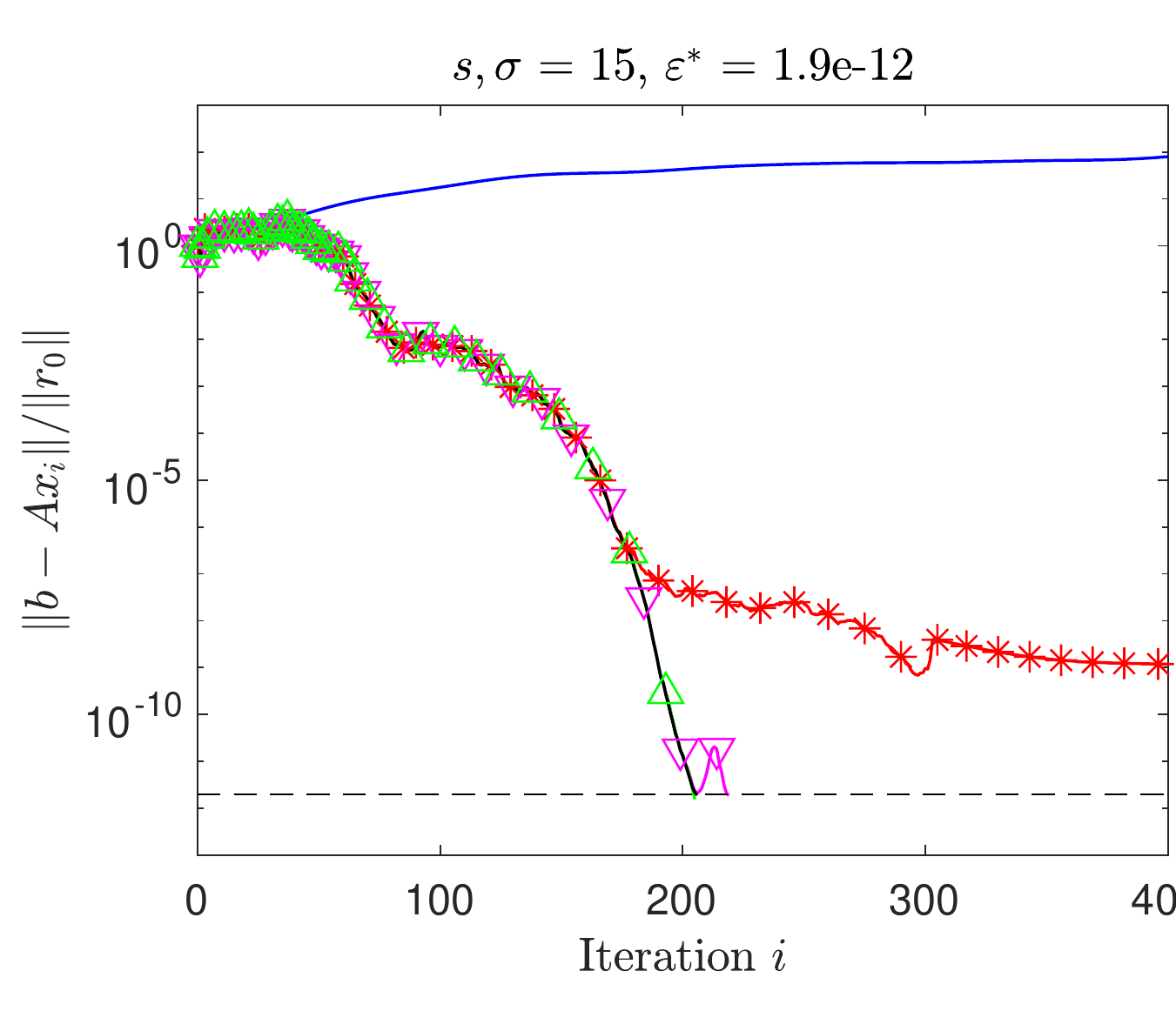}
\includegraphics[width=2.5in,trim=1.25in 3.2in 1.85in 3.15in, clip]
	{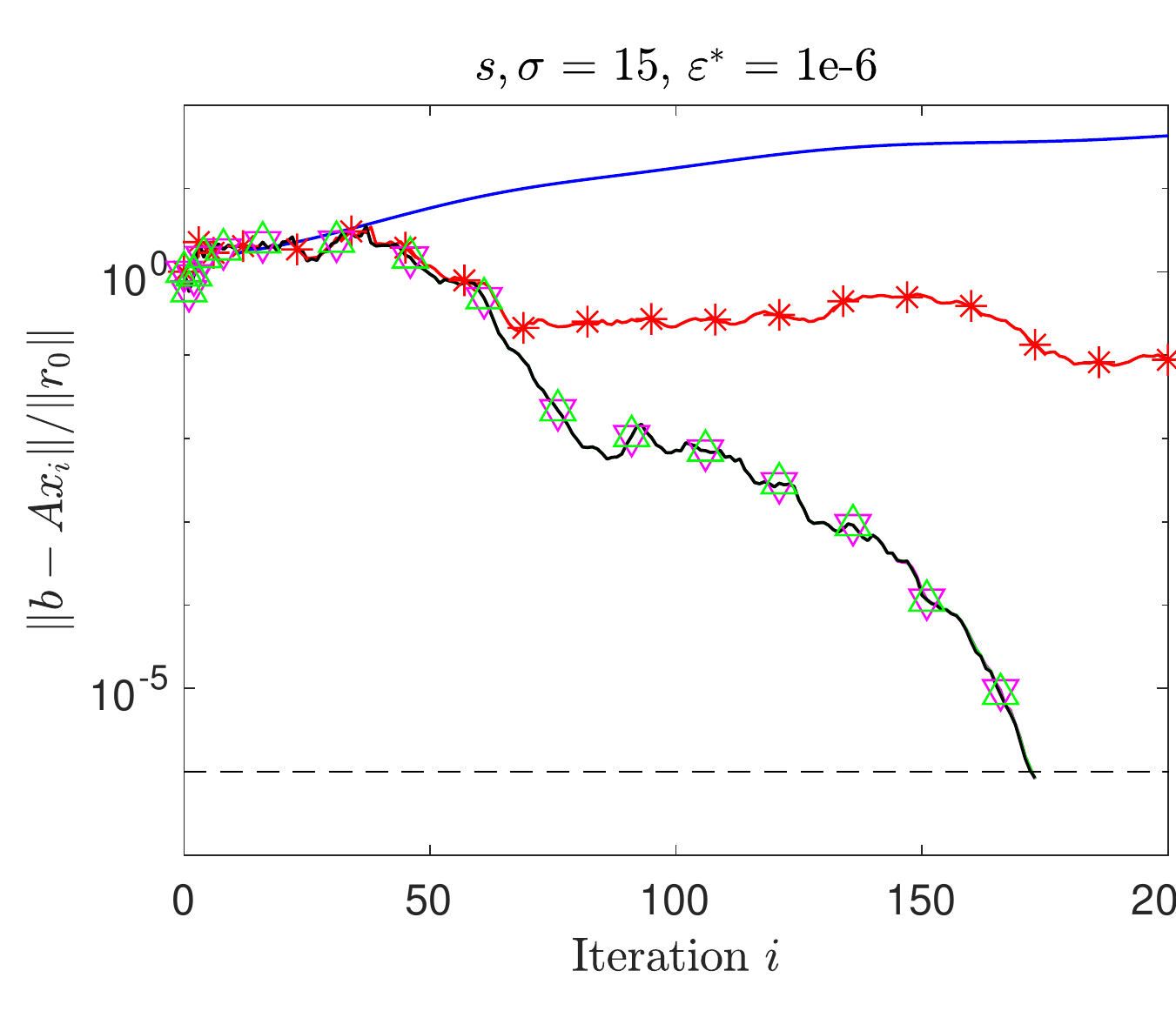} 
\caption{\footnotesize{Convergence of the relative true residual 2-norm for the matrix \texttt{bcsstk09}. 
}}
\label{fig.bcsstk09n_0}
\end{figure}
\begin{table}[hbt!]
\footnotesize
\centering
\resizebox{\textwidth}{!}{%
\begin{tabular}{ll|c|c|c|c|c|}
\cline{3-7}
 &  & fixed s-step CG & old adptv s-step CG & \begin{tabular}[c]{@{}c@{}}impr. adptv s-step\\ CG w/Newton basis\end{tabular} &\begin{tabular}[c]{@{}c@{}}impr. adptv s-step\\ CG w/Chebyshev basis\end{tabular} & HSCG \\ \hline
\multicolumn{1}{|l|}{\multirow{3}{*}{$\varepsilon^*=$ 1.9e-12}} & $s,\sigma=5$ & 51 (253) & 48 (235) & 59 (219) & 54 (205) & \multirow{3}{*}{206} \\ \cline{2-6}
\multicolumn{1}{|l|}{} & $s,\sigma=10$ & - [7e-09] & 48 (375) & 44 (219) & 40 (205) &  \\ \cline{2-6}
\multicolumn{1}{|l|}{} & $s,\sigma=15$ & - & - [2e-10] & 42 (219) & 38 (205) &  \\ \hline
\multicolumn{1}{|l|}{\multirow{3}{*}{$\varepsilon^*=$ 1e-6}} & $s,\sigma=5$ & 35 (173) & 35 (173) & 37 (173) & 37 (173) & \multirow{3}{*}{173} \\ \cline{2-6}
\multicolumn{1}{|l|}{} & $s,\sigma=10$ & 65 (642) & 65 (642) & 21 (173) & 21 (173) &  \\ \cline{2-6}
\multicolumn{1}{|l|}{} & $s,\sigma=15$ & - & 116 (1570) & 16 (173) & 16 (173) &  \\ \hline
\end{tabular}
}

\caption{\footnotesize{Results for experiments with the matrix \texttt{bcsstk09}. 
}}
\label{tab.bcsstk09n0}
\end{table}
}

\afterpage{
\begin{table}[]
\footnotesize
\centering
\resizebox{\textwidth}{!}{%
\begin{tabular}{ll|c|c|c|c|c|}
\cline{3-7}
 &  & fixed s-step CG & old adptv s-step CG & \begin{tabular}[c]{@{}c@{}}impr. adptv s-step\\ CG w/Newton basis\end{tabular} &\begin{tabular}[c]{@{}c@{}}impr. adptv s-step\\ CG w/Chebyshev basis\end{tabular} & HSCG \\ \hline
\multicolumn{1}{|l|}{\multirow{3}{*}{$\varepsilon^*=$ 3.6e-14}} & $s,\sigma=5$ & - [9e-14] & 15 (51) & 23 (51) & 20 (51) & \multirow{3}{*}{52} \\ \cline{2-6}
\multicolumn{1}{|l|}{} & $s,\sigma=10$ & - [1e-10] & - [2e-13] & 21 (51) & 17 (51) &  \\ \cline{2-6}
\multicolumn{1}{|l|}{} & $s,\sigma=15$ & - & - [2e-13] & 21 (51) & 17 (51) &  \\ \hline
\multicolumn{1}{|l|}{\multirow{3}{*}{$\varepsilon^*=$ 1e-6}} & $s,\sigma=5$ & 7 (34) & 7 (34) & 10 (34) & 10 (34) & \multirow{3}{*}{34} \\ \cline{2-6}
\multicolumn{1}{|l|}{} & $s,\sigma=10$ & 5 (50) & 5 (50) & 7 (34) & 7 (34) &  \\ \cline{2-6}
\multicolumn{1}{|l|}{} & $s,\sigma=15$ & - & 19 (263) & 7 (34) & 7 (34) &  \\ \hline
\end{tabular}
}

\caption{\footnotesize{Results for experiments with the matrix \texttt{gr\_30\_30}. 
}}
\label{tab.gr3030n0}
\end{table}

\begin{table}[]
\centering
\resizebox{\textwidth}{!}{%
\begin{tabular}{ll|c|c|c|c|c|}
\cline{3-7}
 &  & fixed s-step CG & old adptv s-step CG & \begin{tabular}[c]{@{}c@{}}impr. adptv s-step\\ CG w/Newton basis\end{tabular} &\begin{tabular}[c]{@{}c@{}}impr. adptv s-step\\ CG w/Chebyshev basis\end{tabular} & HSCG \\ \hline
\multicolumn{1}{|l|}{\multirow{3}{*}{$\varepsilon^*=$ 6.2e-10}} & $s,\sigma=5$ & 32 (158) & 33 (158) & 42 (156) & 33 (119) & \multirow{3}{*}{125} \\ \cline{2-6}
\multicolumn{1}{|l|}{} & $s,\sigma=10$ & - [6e-08] & 54 (455) & 28 (125) & 26 (128) &  \\ \cline{2-6}
\multicolumn{1}{|l|}{} & $s,\sigma=15$ & - & 58 (679) & 29 (156) & 23 (109) &  \\ \hline
\multicolumn{1}{|l|}{\multirow{3}{*}{$\varepsilon^*=$ 1e-6}} & $s,\sigma=5$ & 19 (95) & 19 (95) & 22 (95) & 22 (95) & \multirow{3}{*}{95} \\ \cline{2-6}
\multicolumn{1}{|l|}{} & $s,\sigma=10$ & 52 (511) & 32 (312) & 14 (95) & 13 (95) &  \\ \cline{2-6}
\multicolumn{1}{|l|}{} & $s,\sigma=15$ & - & - & 12 (96) & 11 (96) &  \\ \hline
\end{tabular}
}

\caption{\footnotesize{Results for experiments with the matrix \texttt{nos6}. 
}}
\label{tab.nos6n0}
\end{table}

\begin{table}[]
\centering
\resizebox{\textwidth}{!}{%
\begin{tabular}{ll|c|c|c|c|c|}
\cline{3-7}
 &  & fixed s-step CG & old adptv s-step CG & \begin{tabular}[c]{@{}c@{}}impr. adptv s-step\\ CG w/Newton basis\end{tabular} &\begin{tabular}[c]{@{}c@{}}impr. adptv s-step\\ CG w/Chebyshev basis\end{tabular} & HSCG \\ \hline
\multicolumn{1}{|l|}{\multirow{3}{*}{$\varepsilon^*=$ 2.8e-13}} & $s,\sigma=5$ & - [2e-11] & 67 (217) & 91 (217) & 89 (217) & \multirow{3}{*}{206} \\ \cline{2-6}
\multicolumn{1}{|l|}{} & $s,\sigma=10$ & - [1e-07] & 62 (275) & 80 (217) & 76 (217) &  \\ \cline{2-6}
\multicolumn{1}{|l|}{} & $s,\sigma=15$ & - & - & 78 (217) & 73 (218) &  \\ \hline
\multicolumn{1}{|l|}{\multirow{3}{*}{$\varepsilon^*=$ 1e-6}} & $s,\sigma=5$ & 29 (143) & 29 (143) & 30 (138) & 31 (142) & \multirow{3}{*}{138} \\ \cline{2-6}
\multicolumn{1}{|l|}{} & $s,\sigma=10$ & 31 (309) & 31 (301) & 18 (143) & 18 (139) &  \\ \cline{2-6}
\multicolumn{1}{|l|}{} & $s,\sigma=15$ & - & 443 (5440) & 14 (139) & 14 (140) &  \\ \hline
\end{tabular}
}

\caption{\footnotesize{Results for experiments with the matrix \texttt{mhdb416}. 
}}
\label{tab.mhdb416n0}
\end{table}

\begin{table}[]
\centering
\resizebox{\textwidth}{!}{%
\begin{tabular}{ll|c|c|c|c|c|}
\cline{3-7}
 &  & fixed s-step CG & old adptv s-step CG & \begin{tabular}[c]{@{}c@{}}impr. adptv s-step\\ CG w/Newton basis\end{tabular} &\begin{tabular}[c]{@{}c@{}}impr. adptv s-step\\ CG w/Chebyshev basis\end{tabular} & HSCG \\ \hline
\multicolumn{1}{|l|}{\multirow{3}{*}{$\varepsilon^*=$ 1.1e-09}} & $s,\sigma=5$ & 232 (1159) & 232 (1159) & 221 (1032) & 219 (1035) & \multirow{3}{*}{672} \\ \cline{2-6}
\multicolumn{1}{|l|}{} & $s,\sigma=10$ & - [4e-08] & 547 (5128) & 151 (1047) & 137 (1006) &  \\ \cline{2-6}
\multicolumn{1}{|l|}{} & $s,\sigma=15$ & - & - [2e-07] & 139 (1053) & 123 (1084) &  \\ \hline
\multicolumn{1}{|l|}{\multirow{3}{*}{$\varepsilon^*=$ 1e-6}} & $s,\sigma=5$ & 148 (738) & 148 (738) & 246 (1217) & 219 (1079) & \multirow{3}{*}{509} \\ \cline{2-6}
\multicolumn{1}{|l|}{} & $s,\sigma=10$ & 714 (7134) & 714 (7134) & 134 (1290) & 187 (1820) &  \\ \cline{2-6}
\multicolumn{1}{|l|}{} & $s,\sigma=15$ & - & - & 80 (1086) & 87 (1211) &  \\ \hline
\end{tabular}
}

\caption{\footnotesize{Results for experiments with the matrix \texttt{nos1}. 
}}
\label{tab.nos1n0}
\end{table}
}

\subsection{Discussion}
\label{sec:discussion}

The experimental results clearly show the benefit of the improved adaptive approach over both fixed $s$-step CG and the old adaptive $s$-step CG. 
First, the improved algorithm is more reliable. Whereas a poor choice of $s$ or $\sigma$ can cause failure to converge to the prescribed accuracy in the other 
algorithms, the improved algorithm converges to the desired tolerance in all tested cases. The number of outer loop iterations (global synchronizations) decreases with increased $\sigma$ in the improved approach, and in all cases, this cost is reduced relative to HSCG. Further, the improved algorithm does not require a heuristic choice of $c_{m+j+1}$, but instead automatically sets this parameter based on an inexpensive computation using quantities already available. 
Another benefit of the improved approach is that unlike the old adaptive approach, the total number of iterations required for convergence is independent of $\sigma$.  

In the improved approach, we 
use information generated during the iterations to construct more well-conditioned polynomial bases, which allows for larger $s_k$ values to be used 
while still achieving the desired accuracy. We note that the adaptive approaches (both old and new) are developed \emph{only} with the goal of attaining a prescribed accuracy rather than improving the convergence behavior. However, a side effect of using more well-conditioned bases is that the convergence behavior improves and becomes closer to that of HSCG. This is to be expected based on the theoretical analysis of the $s$-step Lanczos algorithm in~\cite{carson15}, which is currently being investigated in order to extend Greenbaum's classic results on the convergence of HSCG in finite precision~\cite{greenbaum89} to $s$-step variants. 

We also note that it is not clear whether the improved algorithm performs better using Newton or Chebyshev bases. In many cases, we observe that the Chebyshev basis provides the fastest convergence (even faster than HSCG in some cases), however, this is not universal; see, for example, the experiments for \texttt{nos1} with $\sigma=10$ and $\veps^*=10^{-6}$ (the final row of plots in Figure~\ref{fig.nos1test}). 

We note that in most experiments, the number of outer loop iterations required by fixed $s$-step CG, old adaptive $s$-step CG, and improved adaptive $s$-step CG is about the same for $s,\sigma=5$, with the fixed $s$-step and old adaptive approaches outperforming the improved approach in a few cases. This suggests that for small values like $s,\sigma=5$, it is unlikely to be beneficial to use the improved adaptive approach. In this case, a simple monomial basis is likely good enough to ensure that the desired accuracy is attained. Depending on the structure and size of the problem and the particular hardware on which the problem is solved, it may be that a small $s$ value is the best option for minimizing the time per iteration anyway. Providing a concrete answer regarding which approach to use when is nearly impossible due to the highly problem-dependent nature of the performance of iterative methods. Here we have only presented a few experiments on small matrices in order to demonstrate the numerical behavior of the algorithms. In order to provide a clearer picture, we stress that these algorithms should be implemented and compared on a high-performance parallel machine for problems from a wide variety of domains.

\section{Conclusion and future work}
\label{sec:conc}

In this work, we presented an improved adaptive $s$-step CG algorithm for solving SPD linear systems. 
The primary improvement over the previous adaptive $s$-step CG algorithm is due to the use of the method of Meurant and Tich{\' y}~\cite{meti18} for incremental estimation of the 
largest and smallest Ritz values, which are used to dynamically improve the conditioning of the $s$-step basis matrices (which has a large influence on the numerical behavior of the $s$-step methods) as the iterations proceed and also to automatically set a parameter that represents the ratio between the size of the error and the size of the residual. The improved algorithm can provide convergence behavior closer to that of HSCG while also allowing for the use of larger $s_k$ values without unacceptable loss of accuracy. We made an additional small improvement to the criterion for determining $s_k$ based on the observation that only a subset of the basis vectors are needed to compute the iterate updates in each subsequent step.
Our numerical experiments verify that the improved algorithm provides increased reliability and, in many cases, a reduction in the total number of global synchronizations required to converge to the prescribed tolerance versus the old adaptive $s$-step CG algorithm.  

The adaptive approach in $s$-step CG algorithms is an example of how finite precision analysis can provide insight leading to more numerically stable algorithms. It remains to show, however, that the adaptive $s$-step approaches provide a benefit to performance over classical Krylov subspace methods like HSCG. In practical implementations, there remain a few  parameters in the algorithm which must be tuned to maximize performance, including the maximum $s_k$ value $\sigma$ and the maximum basis growth factor $f$. A high-performance parallel implementation and a thorough exploration of the design space remains critical future work.

 \bibliographystyle{siam}
 \bibliography{paper}

%
%
%
%
%
%
%
%
%
%
%
%
%
%
%
%

\end{document}